\pdfoutput=1
\documentclass[aop]{imsart}
\RequirePackage{natbib}
\usepackage{amsmath,amssymb}
\usepackage[mathscr]{eucal}

\usepackage{mathtools}                                                      
\mathtoolsset{showonlyrefs=true}                                    

\usepackage{pdfsync}

\let\counterwithin\relax
\usepackage{chngcntr}

\usepackage{ulem}

\usepackage{amsmath}                                                        
\usepackage{graphicx}                                                       
\usepackage{subfigure}
\usepackage{enumitem}                                                    
\usepackage{hyperref}
\usepackage{amssymb}                                                        
\usepackage[mathscr]{eucal}                                             
\usepackage{cancel}                                                             
\usepackage{pstricks}
\usepackage{rotating}
\usepackage{lscape}
\usepackage[paperwidth=8.5in,paperheight=11in,top=1.25in, bottom=1.25in, left=1.00in, right=1.00in]{geometry}
\usepackage{mathtools}                                                      
\mathtoolsset{showonlyrefs=true}                                    
\usepackage{fixltx2e,amsmath}                                           
\MakeRobust{\eqref}
\usepackage{mathdots}
\usepackage{amsthm}                                                             
\allowdisplaybreaks                                                             
\usepackage{chngcntr}

\usepackage{ulem}       
\usepackage{mathabx}

\usepackage{titlesec}
\usepackage[titletoc,toc,title]{appendix}

\usepackage{apptools, etoolbox}
\usepackage{pst-pdf}

\AtAppendix{\counterwithin{theorem}{section}}

\usepackage{color}

\newtheorem{theorem}{Theorem}[section]
\newtheorem{proposition}[theorem]{Proposition}

\newtheorem{lemma}[theorem]{Lemma}
\newtheorem{definition}[theorem]{Definition}
\newtheorem{example}[theorem]{Example}
\newtheorem{remark}[theorem]{Remark}
\newtheorem{assumption}[theorem]{Assumption}

\def \Qg {Y}

\def \gvar {\gamma}
\def \gY {\gamma^{B}}
\def \gIW {\gamma^{\text{\rm\tiny IW}}}
\def \gIWi {\gamma^{\text{\rm\tiny IW},-1}}
\def \au {\bar{a}}
\def \bu {\bar{b}}

\def \LL {L}

\def \AA {\mathbf{A}}
\def \MM {\mathbf{m}}

\def \o {{\omega}}
\def \a {{\alpha}}
\def \b {{\beta}}
\def \d {{\delta}}
\def \l {{\lambda}}

\def \pG {{\Gamma}}
\def \G {{\Gamma}}
\def \s {{\sigma}}
\def \w {{\omega}}

\def \R {{\mathbb {R}}}
\def \N {{\mathbb {N}}}

\def \x {{\xi}}
\def \e {{\varepsilon}}

\def \r {{\varrho}}

\def \t {{\tau}}
\def \t {{\tau}}

\def \m {{\mu}}
\def \y {{\eta}}

\def \z {{\zeta}}

\def \g {{\gamma}}
\def \O {{\Omega}}
\def \phi {{\varphi}}

\def \tilde {\widetilde}

\def\p{\partial}

\def \B {{\cal{B}}}

\def \o {{\omega}}
\def \a {{\alpha}}
\def \b {{\beta}}
\def \d {{\delta}}
\def \l {{\lambda}}

\def \G {{\Gamma}}
\def \s {{\sigma}}
\def \w {{\omega}}

\def \R  {{\mathbb {R}}}
\def \x {{\xi}}
\def \g {{\gamma}}
\def \e {{\varepsilon}}

\def \t {{\tau}}

\def \m {{\mu}}
\def \y {{\eta}}

\def \z {{\zeta}}
\def \p {{\partial}}
\def \a {{\alpha}}
\def \O {{\Omega}}

\def \d {{\delta}}

\def \o {{\omega}}
\def \a {{\alpha}}
\def \b {{\beta}}
\def \d {{\delta}}

\def \G {\Ga}

\def \Ga {{\Gamma}}

\def \s {{\sigma}}

\def \w {{\omega}}
\def \R {{\mathbb {R}}}
\def \N {{\mathbb {N}}}

\def \x {{\xi}}
\def \e {{\varepsilon}}

\def \r {{\varrho}}

\def \t {{\tau}}
\def \t {{\tau}}

\def \m {{\mu}}
\def \y {{\eta}}

\def \z {{\zeta}}

\def \g {{\gamma}}
\def \O {{\Omega}}
\def \phi {{\varphi}}

\def \tilde {\widetilde}

\def\l {\lambda}

\def \F {\mathcal{F}}
\def \B {\mathscr{B}}

\def \à {{\`a }}
\def \è {{\`e }}
\def \ò {{\`o }}
\def \ù {{\`u }}

\begin{document}

%
%
%
%
%

\begin{frontmatter}

\title{On stochastic Langevin and Fokker-Planck equations: the two-dimensional case}
\runtitle{On stochastic Langevin and Fokker-Planck equations}

\begin{aug}
  \author{\fnms{Andrea}  \snm{Pascucci}\corref{}
  \ead[label=e1]{andrea.pascucci@unibo.it}}
  \and
  \author{\fnms{Antonello}  \snm{Pesce}%
  \ead[label=e2]{antonello.pesce2@unibo.it}%
  \ead[label=u1,url]{http://www.foo.com}}


  \runauthor{A. Pascucci et al.}

  \affiliation{Universit\`a di Bologna}

  \address{Dipartimento di Matematica, Piazza di Porta san Donato, 5 Bologna (Italy),\\
          \printead{e1,e2}}

\end{aug}

\begin{abstract}
We prove existence, regularity in H\"older classes and estimates from above and below of the
fundamental solution of the stochastic Langevin equation. This degenerate SPDE satisfies the weak
H\"orman\-der condition. We use a Wentzell's transform to reduce the SPDE to a PDE with random
coefficients; then we apply a new method, based on the parametrix technique, to construct a
fundamental solution. This approach avoids the use of the Duhamel's principle for the SPDE and the
related measurability issues that appear in the stochastic framework. Our results are new even for
the deterministic equation.
\end{abstract}



\end{frontmatter}

%
%

\section{Introduction}\label{intro}

%
%
%
%

We consider the stochastic version of the Fokker-Planck equation
\begin{equation}\label{PDE1}
 \p_{t}u+\sum_{j=1}^{n}v_{j}\p_{x_{j}}u=\frac{1}{2}\sum_{i,j=1}^{n}a_{ij}\partial_{v_{i}v_j}u.
\end{equation}
Here the variables $t\ge 0$, $x\in\R^{n}$ and $v\in\R^{n}$ respectively stand for time, position
and velocity, and the unknown $u=u_{t}(x,v)\ge 0$ stands for the density of particles in phase
space. The vector field
  $\mathbf{Y}:=\p_{t}+v\cdot\nabla_{x}$
on the left-hand side of \eqref{PDE1} describes transport; the coefficients
$a_{ij}$ describe some kind of collision among particles and in general may depend on the solution
$u$ through some integral expressions. Linear Fokker-Planck equations (cf. \cite{Desvillettes} and
\cite{Risken}), non-linear Boltzmann-Landau equations (cf. \cite{Lions1} and \cite{Cercignani})
and non-linear equations for Lagrangian stochastic models commonly used in the simulation of
turbulent flows (cf. \cite{Talay}) can be written in the form \eqref{PDE1}. In mathematical
finance, \eqref{PDE1} describes path-dependent financial contracts such as Asian options (see, for
instance, \cite{Pascucci2011}).

In this note we study a kinetic model where the position and the velocity of a particle are
stochastic processes $(X_{t},V_{t})$ only partially observable through some observation process
$O_{t}$. {We consider the two-dimensional case, $n=1$, which is already challenging enough, and
propose an approach that hopefully can be extended to the multi-dimensional case}. If
$\F^{O}_{t}=\s(O_{s},\, s\le t)$
denotes the filtration of the observations then, under natural assumptions, 
the {\it conditional} density  $p_{t}(x,v)$ of $(X_{t},V_{t})$ given $\F_{t}^{O}$ solves a linear
SPDE of the form
\begin{align}\label{spde}
  d_{\mathbf{Y}}u_{t}(x,v)=\frac{a_t(x,v)}{2}\p_{vv}u_t(x,v)dt+\s_t(x,v)\p_vu_t(x,v)dW_t,\qquad
  \mathbf{Y}=\p_{t}+v\p_{x}.
\end{align}
In \eqref{spde} $W$ is a Wiener process defined on a complete probability space $(\O,\F,P)$
endowed with a filtration $\left(\F_{t}\right)_{t\ge 0}$ satisfying the usual conditions. The
symbol $d_{\mathbf{Y}}$ indicates that the equation is solved in the It\^o (or strong) sense: a
solution to \eqref{spde} is a continuous process
$u_{t}=u_{t}(x,v)$ that is twice differentiable in $v$ and 
such that
\begin{align} \label{spde1}
 u_{t}\left(\gY_{t}(x,v)\right)=u_{0}(x,v)+\frac{1}{2}\int_{0}^{t}\left(a_{s}\p_{vv}u_{s}\right)
 (\gY_{s}(x,v))ds+
 \int_{0}^{t}\left(\s_{s}\p_vu_{s}\right)(\gY_{s}(x,v))dW_{s}
\end{align}
where $t\mapsto\gY_{t}(x,v)$ denotes the integral curve, starting from $(x,v)$, of the advection
vector field $v\p_{x}$, that is 
\begin{equation}\label{e4}
 \gY_{t}(x,v)=e^{t B}(x,v)=(x+tv,v),\qquad B=\begin{pmatrix}
 0 & 1 \\ 0 & 0
 \end{pmatrix}.
\end{equation}
Clearly, in case the observation process $O$ is independent of $X$ and $V$, the SPDE \eqref{spde}
boils down to the deterministic PDE \eqref{PDE1} with $n=1$.

The main goal of this paper is to show existence, regularity and Gaussian-type estimates of a
stochastic fundamental solution of \eqref{spde}. As far as the authors are aware, such kind of
results was never proved for SPDEs that satisfy the {\it weak} H\"ormander condition, that is
under the assumption that the drift has a key role in the noise propagation. We mention that
hypoellipticity for SPDEs under the strong H\"ormander condition was studied by Chaleyat-Maurel
and Michel \cite{MR736147}, Kunita \cite{MR705933}, Krylov \cite{Krylov17} and Jinniao
\cite{MR3706782}. Even in the deterministic case, our results are new in that they extend the
recent results \cite{MR2659772}, \cite{MR3758337} for Kolmogorov equations with general drift.

Our method is based on a Wentzell's reduction of the SPDE to a PDE with random coefficients to
which we apply the parametrix technique to construct a fundamental solution. This approach avoids
the use of the Duhamel's principle for the SPDE and the related measurability issues that appear
in the stochastic framework as discussed, for instance, in \cite{Sowers}. As in
\cite{PascucciPesce1}, Wentzell's reduction of the SPDE is done globally: to control the behavior
as $|x|,|v|\rightarrow\infty$ of the random coefficients of the resulting PDE, we impose some
flattening condition at infinity on the coefficient $\s_{t}(x,v)$ in \eqref{spde} (cf. Assumption
\ref{Ass3}). Compared to the uniformly parabolic case, two main new difficulties arise:
\begin{itemize}
  \item[i)] the It\^o-Wentzell transform drastically affects the drift $\mathbf{Y}$: in particular,
  after the random change of coordinates, the new
  drift has no longer polynomials coefficients. Consequently, a careful analysis is needed to
  check the validity of the H\"ormander condition in the new coordinates.
  This question is discussed in more detail in Section \ref{sub01};
  \item[ii)] in the deterministic case, the parametrix method
  has been applied to degenerate Fokker-Planck equations,
  including \eqref{spde} with $\s\equiv 0$, by several authors, \cite{Polidoro2}, \cite{DiFrancescoPascucci2}, \cite{MR2802040},
  \cite{MR3891754},
  using {\it intrinsic} H\"older spaces. Loosely speaking, the intrinsic H\"older regularity reflects the
  geometry of the PDE and is defined in terms of the translations and homogeneous norm associated to
  the H\"ormander vector fields: this kind of regularity is natural for the study
  of the singular kernels that come into play in the parametrix iterative procedure. Now, under the weak H\"ormander condition,
  the intrinsic regularity properties in space and time are closely intertwined and cannot be studied
  separately.
  However, assuming that the coefficients are
  merely predictable, we have no good control on the regularity in the time variable; for instance,
  even in the deterministic case, the coefficients are only measurable in $t$ and consequently they cannot be H\"older continuous in $(x,v)$ in the intrinsic sense.
  On the other hand, assuming that the coefficients are H\"older continuous in $(x,v)$ in the
  classical Euclidean sense, the parametrix method still works as long as we use a suitable
  {\it time-dependent} parametrix and exploit the fact
  that the intrinsic translations coincide with the Euclidean ones
  for points $(t,x,v)$ and $(t,\x,\y)$ at the same time level.
  We comment on this question more thoroughly in Section \ref{sub02}.
\end{itemize}
The rest of the paper is organized as follows. In Sections \ref{sub01} and \ref{sub02} we go
deeper into the issues mentioned above. In Section \ref{results} we set the assumptions, introduce
the functional setting and state the main result, Theorem \ref{t1}. In Section \ref{pointw} we
prove some crucial estimate for stochastic flows of diffeomorphisms: these estimates, which can be
of independent interest, extend some result of \cite{MR1070361}. In Section \ref{ItoWen} we
formulate a version of the It\^o-Wentzell formula and exploit the results of Section \ref{pointw}
to perform a stochastic change of variable in order to reduce the SPDE to a PDE with random
coefficients. In Section \ref{mproof} we build on the work by Delarue and Menozzi \cite{MR2659772}
to develop a parametrix method for Kolmogorov PDEs with general drift (Theorem \ref{t3}). Finally,
in Section \ref{proof} we complete the proof of Theorem \ref{t1}.

\subsection{Stochastic Langevin equation and the H\"ormander condition}\label{sub01}
For illustrative purposes, we examine the case of constant coefficients and introduce the
stochastic counterpart of the classical Langevin PDE.

Let $B,W$ be independent real Brownian motions, $a>0$ and $\s\in[0,\sqrt{a}]$. The Langevin model
is defined in terms of the system of SDEs
\begin{equation}\label{SDELang}
\begin{cases}
    dX_{t}=V_{t}dt,\\
    dV_{t}=\sqrt{a-\s^2} dB_{t} - \s dW_{t}.
  \end{cases}
\end{equation}
We interpret $W$ as the observation process: if $\s=0$ the velocity $V$ is unobservable, while for
$\s=\sqrt{a}$ the velocity $V$ is completely observable, being equal to $W$. To shorten notations,
we denote by $z=(x,v)$ and $\z=(\x,\y)$ the points in $\R^{2}$. Setting $Z_{t}=(X_{t},Y_{t})$,
equation \eqref{SDELang} can be rewritten as
\begin{equation}\label{SDELang1}
 dZ_{t}=BZ_{t}dt+\mathbf{e}_{2} d(\sqrt{a-\s^2} B_{t} - \s W_{t}),\qquad \mathbf{e}_{2}={0\choose 1},
\end{equation}
where $B$ is the matrix in \eqref{e4}.

In this section we show in two different ways that the SPDE
\begin{equation}\label{e2}
  d_{\mathbf{Y}}u_{t}=
  \frac{a}{2}\p_{vv}u_{t}dt+\s \p_{v}u_{t}dW_{t},\qquad \mathbf{Y}:=\p_{t}+v\p_{x},
\end{equation}
is the {\it forward Kolmogorov} (or {\it Fokker-Planck}) equation of the SDE \eqref{SDELang}
conditioned to the Brownian observation given by $\F^{W}_{t}=\s(W_{s},\, s\le t)$. In the
uniformly parabolic case, this is a well-known fact, proved under diverse assumptions by several
authors (see, for instance, \cite{MR0242552}, \cite{MR0501350} and \cite{MR553909}).

In the first approach, we solve explicitly the linear SDE \eqref{SDELang1} and find the expression
of the conditional transition density $\pG$ of the solution $Z$: by It\^o formula, we directly
infer that $\pG$ is the fundamental solution of the SPDE \eqref{e2}. The second approach, inspired
by \cite{MR1795614}, is much more general because it does not require the explicit knowledge of
$\pG$: we first prove the existence of the fundamental solution of the SPDE \eqref{e2} and then
show that it is the conditional transition density of the solution of \eqref{SDELang}.

The following result is an easy consequence of the It\^o formula and isometry.
\begin{proposition}
The solution $Z=Z^{\z}$ of \eqref{SDELang1}, with initial condition $\z=(\x,\y)\in\R^{2}$, is
given by
  $$Z^{\z}_{t}=e^{tB}\left(\z+\int_{0}^{t}e^{-sB}\mathbf{e}_{2}\, d(\sqrt{a-\s^2} B_{s} - \s W_{s})\right)$$
with $\mathbf{e}_{2}$ as in \eqref{SDELang1}. Conditioned to $\F_{t}^{W}$, $Z^{\z}_{t}$ has normal
distribution with mean and covariance matrix given by
\begin{align}\label{e7b}
   m_{t}(\z)&:=E\left[Z^{\z}_{t}\mid\F_{t}^{W}\right]=e^{tB}\left(\z-\s\int_{0}^{t}e^{-sB}\mathbf{e}_{2} dW_{s}\right)=
   \begin{pmatrix}
     \x+t\y-\s\int_{0}^{t}(t-s)dW_{s} \\
     \y-\s W_{t} \
   \end{pmatrix},\\ \label{e7c}
   \mathcal{C}_{t}&:=\text{cov}\left(Z^{\z}_{t}\mid\F_{t}^{W}\right)=(a-\s^{2})Q_{t},\qquad Q_{t}:=
   \int_{0}^{t}\left(e^{sB}\mathbf{e}_{2}\right)
   \left(e^{sB}\mathbf{e}_{2}\right)^{\ast}ds=\begin{pmatrix}
    \frac{t^{3}}{3} & \frac{t^{2}}{2} \\
    \frac{t^{2}}{2} & t \
  \end{pmatrix}.
\end{align}
In particular, if $\s=\sqrt{a}$ then the distribution of $Z^{\z}_{t}$ conditioned to $\F_{t}^{W}$
is a Dirac delta centered at $m_{t}(\z)$; if $\s\in[0,\sqrt{a})$ and $t>0$ then $Z^{\z}_{t}$ has
density, conditioned to $\F_{t}^{W}$, given by
\begin{align}\label{e7}
 \pG(t,z;0,\z)&=\frac{1}{2\pi\sqrt{\det\mathcal{C}_{t}}}\exp\left(-\frac{1}{2}\langle \mathcal{C}_{t}^{-1}(z-m_{t}(\z)),
 (z-m_{t}(\z))\rangle\right),\qquad z\in\R^{2}.
\end{align}
More explicitly, we have $\pG(t,z;0,\z)=\pG_{0}(t,z-m_{t}(\z))$ where
\begin{equation}\label{e6}
   \pG_{0}(t,x,v)=\frac{\sqrt{3}}{\pi t^{2}(a-\s^{2})}\exp\left(-\frac{2}{a-\s^{2}}\left(\frac{v^{2}}{t}-
 \frac{3vx}{t^{2}}+\frac{3x^{2}}{t^{3}}\right)\right),\qquad t>0,\ (x,v)\in\R^{2}.
\end{equation}
By the It\^o formula, $\pG(t,z;0,\z)$ is the stochastic fundamental solution of SPDE \eqref{e2},
with pole at $(0,\z)$.
\end{proposition}

As an alternative approach, we construct the fundamental solution of the SPDE \eqref{e2} by
performing some suitable change of variables. First we transform \eqref{e2} into a PDE with random
coefficients, satisfying the weak H\"ormander condition; {by a second change of variables, we
remove the drift of the equation and
transform it into a deterministic heat equation. 
Going back to the original variables, we find the stochastic fundamental solution of \eqref{e2},
which obviously coincides with $\pG$ in \eqref{e7}. Eventually, we prove that $\pG(t,\cdot;0,\z)$
is a density of $Z^{\z}_{t}$ conditioned to $\F^{W}_{t}$. We split the proof in three steps.}

\smallskip
\noindent{\bf [Step 1]}\ We set
\begin{equation}\label{e15}
  \hat{u}_{t}(x,v)=u_{t}(x,v-\s W_{t}).
\end{equation}
By It\^o formula, $u$ solves \eqref{e2} if and only if $\hat{u}$ solves the Langevin PDE
 \begin{equation}\label{e3bis}
 \p_{t}\hat{u}+(v-\s W_{t})\p_{x}\hat{u}=\frac{a-\s^{2}}{2}\p_{vv}\hat{u}.
 \end{equation}
By this change of coordinates we get rid of the stochastic part of the SPDE; however, this is done
at the cost of introducing a random drift term. For the moment, this is not a big issue because
$\s$ is constant and, in particular, independent of $v$: for this reason, the weak H\"ormander
condition is preserved since the vector fields $\p_{v}$, $\p_{t}+(v-\s W_{t})\p_{x}$ and their Lie
bracket
  $$[\p_{v},\p_{t}+(v-\s W_{t})\p_{x}]=\p_{x}$$
span $\R^{3}$ at any point.

\smallskip \noindent{\bf [Step 2]}\ In order to remove the random drift, we perform a second change
of variables:
\begin{equation}\label{e39}
  g_{t}(x,v)=\hat{u}_{t}(\gvar_t(x,v)),\qquad \gvar_{t}(x,v):=\left(x+tv-\s\int_{0}^{t}W_{s}ds,v\right).
\end{equation}
The spatial Jacobian of $\gvar_{t}$ equals
 $$\nabla \gvar_{t}(x,y)=\begin{pmatrix}
  1 & t \\
  0 & 1
 \end{pmatrix}$$
so that $\gvar_{t}$ is one-to-one and onto $\R^{2}$ for any $t$. Then, \eqref{e3bis} is
transformed into the deterministic heat equation with time-dependent coefficients
\begin{equation}\label{e3b}
 \p_{t}g_{t}(x,v)=\frac{a-\s^{2}}{2}\left(t^{2}\p_{xx}-2t\p_{xv}+\p_{vv}\right)g_{t}(x,v).
\end{equation}
Equation \eqref{e3b} is not uniformly parabolic because the matrix of coefficients of the second
order part
  $$a_{t}:=(a-\s^{2})\begin{pmatrix}
    t^2 & -t \\
    -t & 1 \
  \end{pmatrix}$$
is singular. However, in case of partial observation, that is $\s\in[0,\sqrt{a})$, the diffusion
matrix
\begin{equation}\label{e32}
  A_{t}=\int_{0}^{t}a_{s}ds=(a-\s^{2})\begin{pmatrix}
    \frac{t^{3}}{3} & -\frac{t^{2}}{2} \\
    -\frac{t^{2}}{2} & t \
  \end{pmatrix}
\end{equation}
is positive definite for any $t>0$ and therefore \eqref{e3b} admits a Gaussian fundamental solution.
For $\s=0$, this result was originally proved by Kolmogorov \cite{Kolmogorov2} (see also the
introduction in \cite{Hormander}). Going back to the original variables we recover the explicit
expression of $\pG$ in \eqref{e7}.

Incidentally, we notice that \eqref{e3b} also reads
\begin{equation}\label{e3c}
  \p_{t}g_{t}(x,v)=\frac{a-\s^{2}}{2}\bar{\mathbf{V}}_{t}^{2}g_{t}(x,v),\qquad \bar{\mathbf{V}}_{t}:=\p_{v}-t\p_{x},
\end{equation}
where the vector fields $\p_{t}$ and $\bar{\mathbf{V}}_{t}$ satisfy the weak H\"ormander condition
in $\R^{3}$ because $[\bar{\mathbf{V}}_{t},\p_{t}]=\p_{x}$.\medskip

\noindent{\bf [Step 3]}\ We show that $\pG$ is the conditional transition density of $Z$: the
proof is based on a combination of the arguments of \cite{MR1795614} with the gradient estimates
for Kolmogorov equations proved in \cite{MR2352998}.
\begin{theorem}
Let $Z^{\z}$ denote the solution of the linear SDE \eqref{SDELang1} starting from $\z\in\R^{2}$
and let $\pG=\pG(t,\cdot;0,\z)$ in \eqref{e7} be the fundamental solution of the Langevin SPDE
\eqref{e2} with $\s\in[0,\sqrt{a})$. For any bounded and measurable function $\phi$ on $\R^{2}$,
we have
\begin{equation}\label{e5}
 E\left[\phi(Z^{\z}_{t})\mid \F^W_t\right]=\int_{\R^2}\phi(z)\pG(t,z;0,\z)dz.
\end{equation}
\end{theorem}
\proof It is not restrictive to assume that $\phi$ is a test function. Let
  $$I_{t}(\z):=\int_{\R^2}\phi(z)\pG(t,z;0,\z)dz,\qquad t>0,\ \z\in\R^{2}.$$
By \eqref{e7b}-\eqref{e7}, $I_{t}(\z)$ is $\F^{W}_{t}$-measurable: thus, to prove the thesis it
suffices to show that, for any continuous and non-negative function $c=c_{s}(w)$ on
$[0,t]\times\R$, we have
\begin{equation}\label{e5b}
  E\left[e^{-\int_{0}^{t}c_s(W_s)ds}\phi(Z^{\z}_t)\right]=
  E\left[e^{-\int_{0}^{t}c_s(W_s)ds}I_{t}(\z)\right].
\end{equation}
Let
\begin{equation}\label{e8}
  \mathcal{L}^{(\s)}=\frac{a}{2}\left(\p_{vv}-2\s\p_{vw}+\p_{ww}\right)+v\p_x
\end{equation}
be the infinitesimal generator of the three-dimensional process $(X,V,W)$. For
$\s\in[0,\sqrt{a})$, $\p_{t}+\mathcal{L}^{(\s)}$ satisfies the weak H\"ormander condition in
$\R^{4}$ and has a Gaussian fundamental solution (see, for instance, formula (2.9) in
\cite{MR2352998}). We denote by $f=f_{s}(x,v,w)$ the classical solution of the {\it backward}
Cauchy problem
\begin{equation}\label{e9}
  \begin{cases}
    \left(\p_{s}+\mathcal{L}^{(\s)}\right)f_{s}(x,v,w)-c_{s}(w)f_{s}(x,v,w)=0, & (s,x,v,w)\in[0,t)\times\R^{3}, \\
    f_{t}(x,v,w)=\phi(x,v), & (x,v,w)\in\R^{3},
  \end{cases}
\end{equation}
and set
  $$M_{s}:=e^{-\int_{0}^{s} c_{\t}(W_{\t})d\t}\int_{\R^2}f_s(z,W_s)\pG(s,z;0,\z)dz,\qquad s\in [0,t].$$
By definition, we have
\begin{align}
 E\left[M_{t}\right]&=E\left[e^{-\int_{0}^{t}c_s(W_s)ds}I_{t}(\z)\right]
\intertext{and, by the Feynman-Kac representation of $f$,}
 E\left[M_{0}\right]&=f_0(\z,0)=E\left[e^{-\int_{0}^{t} c_{s}(W_{s})d s}\phi(Z^{\z}_{t})\right].
\end{align}
Hence \eqref{e5b} follows by proving that $M$ is a martingale. By the It\^o formula, we have
\begin{align}
 df_{s}(x,v,W_{s})&= \left(\p_{s}f_{s}+\frac{1}{2}\p_{ww}f_{s}\right)(x,v,W_{s})ds+
 \left(\p_{w}f_{s}\right)(x,v,W_{s})dW_{s}\\
 &= \left(-\mathcal{L}^{(\s)}f_{s}+c_{s}f_{s}+\frac{1}{2}\p_{ww}f_{s}\right)(x,v,W_{s})ds+
 \left(\p_{w}f_{s}\right)(x,v,W_{s})dW_{s}.
\end{align}
Moreover, since $\pG$ solves the SPDE \eqref{e2}, setting $e_{s}:=e^{-\int_{0}^{s}
c_{\t}(W_{\t})d\t}$ for brevity, we get
\begin{align}
 dM_{s}=-c_{s}(W_{s})M_{s}ds&+e_{s}\int_{\R^{2}}\left(-\mathcal{L}^{(\s)}f_{s}+c_{s}f_{s}
 +\frac{1}{2}\p_{ww}f_{s}\right)(x,v,W_{s})\pG(s,x,v;0,\z)dxdv\,ds\\
 &+e_{s}\int_{\R^{2}}\left(\p_{w}f_{s}\right)(x,v,W_{s})\pG(s,x,v;0,\z)dxdv\,dW_{s}\\
 &+e_{s}\int_{\R^{2}}f_{s}(x,v,W_{s})\left(\frac{a}{2}\p_{vv}-v\p_{x}\right)\pG(s,x,v;0,\z)dxdv\,ds\\
 &+e_{s}\s\int_{\R^{2}}f_{s}(x,v,W_{s})\p_{v}\pG(s,x,v;0,\z)dxdv\,dW_{s}\\
 &+e_{s}\s\int_{\R^{2}}\p_{w}f_{s}(x,v,W_{s})\p_{v}\pG(s,x,v;0,\z)dxdv\,ds.
\end{align}
Integrating by parts, we find
  $$dM_{s}=e_{s}\int_{\R^{2}}
  \left(\p_{w}f_{s}-\s\p_{v}f_{s}\right)(x,v,W_{s})\pG(s,x,v;0,\z)dxdv\,dW_{s},
  $$
which shows that $M$ is at least a local martingale.

To conclude, we recall the gradient estimates proved in \cite{MR2352998}, Proposition 3.3: for any
test function $\phi$ there exist two positive constants $\e,C$ such that
\begin{equation}\label{e10}
 |\p_{v}f_{s}(x,v,w)|+|\p_{w}f_{s}(x,v,w)|\le \frac{C}{(t-s)^{\frac{1}{2}-\e}},\qquad
 (s,x,v,w)\in[0,t)\times\R^{3}.
\end{equation}
Thus, we have
\begin{align}
  &E\left[\int_{0}^{t}\left(\int_{\R^{2}}
  \left(\p_{w}f_{s}-\s\p_{v}f_{s}\right)(x,v,W_{s})\pG(s,x,v;0,\z)dxdv\right)^{2}ds\right]\\
  &\le \int_{0}^{t}\frac{C}{(t-s)^{1-2\e}}E\left[\left(\int_{\R^{2}}
  \pG(s,x,v;0,\z)dxdv\right)^{2}\right]ds\\
  &=\int_{0}^{t}\frac{C}{(t-s)^{1-2\e}}ds <\infty
\end{align}
and this proves that $M$ is a true martingale.\endproof

\subsection{Intrinsic vs Euclidean H\"older spaces for the deterministic Langevin equation}\label{sub02}
The parametrix method requires some assumption on the regularity of the coefficients of the PDE:
in the uniformly parabolic case, it suffices to assume that the coefficients are bounded, H\"older
continuous in the spatial variables and measurable in time (cf. \cite{Friedman}).

In this paper, we apply the parametrix method assuming that the coefficients of the Langevin SPDE
\eqref{spde} are predictable processes that are H\"older continuous in $(x,v)$ in the Euclidean
sense.
From the analytical perspective this is not the natural choice: indeed, it is
well known that the natural framework for the study 
of H\"ormander operators is the analysis on Lie groups (see, for instance, \cite{MR657581}). In
this section, we motivate our choice to use Euclidean H\"older spaces rather than the intrinsic
ones.

We recall that Lanconelli and Polidoro \cite{LanconelliPolidoro} first studied the intrinsic
geometry of the Langevin operator in \eqref{e2} with $\s=0$:
  $$\mathcal{L}_{a}:=\frac{a}{2}\p_{vv}-v\p_{x}-\p_{t}.$$
They noticed that $\mathcal{L}_{a}$ is invariant with respect to the homogeneous Lie group
$(\R^{3},\ast,\d)$ where the group law is given by
\begin{equation}\label{egl}
  (\t,\x,\y)\ast(t,x,v)=(t+\t,x+\x+t\y,v+\y),\qquad (\t,\x,\y),(t,x,v)\in\R^{3},
\end{equation}
and $\d=(\d_{\l})_{\l>0}$ is the ultra-parabolic dilation operator defined as
\begin{equation}\label{egdil}
  \d_{\l}(t,x,v)=(\l^{2}t,\l^{3}x,\l v),\qquad (t,x,v)\in\R^{3},\ \l>0.
\end{equation}
More precisely, $\mathcal{L}_{a}$ is invariant with respect to the left-$\ast$-translations
$\ell_{(\t,\x,\y)}(t,x,v)=(\t,\x,\y)\ast(t,x,v)$, in the sense that
\begin{equation}
  \mathcal{L}_{a}(f\circ \ell_{(\t,\x,\y)})=\left(\mathcal{L}_{a}f\right)\circ \ell_{(\t,\x,\y)},\qquad
  (\t,\x,\y)\in\R^{3},
\end{equation}
and is $\d$-homogeneous of degree two, in that
\begin{equation}
  \mathcal{L}_{a}(f\circ \d_{\l})=\l^{2}\left(\mathcal{L}_{a}f\right)\circ \d_{\l},\qquad \l>0.
\end{equation}

It is natural to endow $(\R^{3},\ast,\d)$ with the $\d$-homogeneous norm
\begin{equation}\label{ee37}
 |(t,x,v)|_{\mathcal{L}}=|t|^{\frac{1}{2}}+|x|^{\frac{1}{3}}+|v|
\end{equation}
and the distance
\begin{equation}\label{ee37bis}
  d_{\mathcal{L}}\left((t,x,v),(\t,\x,\y)\right)=|(\t,\x,\y)^{-1}\ast(t,x,v)|_{\mathcal{L}}.
\end{equation}
The intrinsic H\"older spaces associated to $d_{\mathcal{L}}$ are particularly beneficial for the
study of existence and regularity properties of solutions to the Langevin equation because they
comply with the asymptotic properties of its fundamental solution $\pG$ near the pole: let us
recall that
 $$\pG(t,x,v;\t,\x,\y)=\pG_{0}\left((\t,\x,\y)^{-1}\ast(t,x,v)\right),\qquad \t<t,$$
where $\pG_{0}$ is the fundamental solution of $\mathcal{L}$ in \eqref{e6} with $\s=0$ and
  $(\t,\x,\y)^{-1}=\left(-\t,-\x+\t\y,-\y\right)$
is the $\ast$-inverse of $(\t,\x,\y)$. Notice also that $\pG$ is $\d$-homogeneous of degree four,
where four is the so-called $\d$-homogeneous dimension of $\R^{2}$.

Based on the use of intrinsic H\"older spaces defined in terms of $d_{\mathcal{L}}$, a stream of
literature has built a complete theory of existence and regularity, analogous to that for
uniformly parabolic PDEs: we mention some of the main contributions \cite{Polidoro2},
\cite{Polidoro1}, \cite{Manfredini}, \cite{Lunardi}, \cite{DiFrancescoPascucci2},
\cite{DiFrancescoPolidoro}, \cite{MR3429628} and, in particular, \cite{Polidoro2},
\cite{DiFrancescoPascucci2}, \cite{Menozzi10} where the parametrix method for Kolmogorov PDEs was
developed.

On the other hand, intrinsic H\"older regularity 
can be a rather restrictive property as shown by the following example.
\begin{example}
For $x,\x\in\R$ and $t\neq \t$, let
\begin{equation}\label{e28}
  z=\left(x,-\frac{x-\x}{t-\t}\right),\qquad \z=\left(\x,-\frac{x-\x}{t-\t}\right)
\end{equation}
Then we have
  $$(\t,\z)^{-1}\ast (t,z)=\left(t-\t,0,0\right),$$
and therefore
  $$d_{\mathcal{L}}((t,z),(\t,\z))=|t-\t|^{\frac{1}{2}}.$$
Since $x$ and $\x$ are arbitrary real numbers, we see that points in $\R^{3}$ that are far from
each other in the Euclidean sense, can be very close in the intrinsic sense. It follows that, if a
function $f(t,x,v)=f(x)$ depends only on $x$ and is H\"older continuous in the intrinsic sense
(i.e. with respect to $d_{\mathcal{L}}$), then it must be constant: in fact, for $z,\z$ as in
\eqref{e28}, we have
  $$\left|f(x)-f(\x)\right|=\left|f(t,z)-f(\t,\z)\right|\le C |t-\t|^{\a}$$
for some positive constants $C,\a$ and for any $x,\x\in\R$ and $t\neq \t$.
\end{example}

When it comes to studying the stochastic Langevin equation, the use of Euclidean H\"older spaces
seems unavoidable. The problem is that we have to deal with functions $f=f_{t}(x,v)$ that are:
\begin{itemize}
  \item[-] {\it H\"older continuous} with respect to the space variables $(x,v)$ in order to apply the parametrix method;
  \item[-] {\it measurable} with respect to the time variable $t$ because $f$ plays the role of a coefficient of the SPDE that is a predictable process
  (i.e. merely measurable in $t$).
\end{itemize}
As opposed to the standard parabolic case, in terms of the metric $d_{\mathcal{L}}$ there doesn't
seem to be a clear way to separate regularity in $(x,v)$ from regularity in $t$: indeed this is
due to the definition of $\ast$-translation that mixes up spatial and time variables (see
\eqref{egl}). On the other hand, we may observe that the Euclidean- and $\ast$- differences of
points at the same time level coincide:
  $$(t,\x,\y)^{-1}\ast(t,x,v)=(0,x-\x,v-\y),\qquad x,v,\x,\y\in\R.$$
Thus, to avoid using $\ast$-translations, the idea is to combine this property with a suitable
definition of {\it time-dependent parametrix} that makes the parametrix procedure work: this will
be done in Section \ref{mproof}.

Concerning the use of the Euclidean or homogeneous norm in $\R^{2}$, let us denote by
$bC^{\a}(\R^{2})$ and $bC^{\a}_{\mathcal{L}}(\R^{2})$ the space of bounded and H\"older continuous
functions with respect to the Euclidean norm and the homogeneous norm $|x|^{\frac{1}{3}}+|v|$,
respectively. Since $|(x,v)|\le |x|^{\frac{1}{3}}+|v|$ for $|(x,v)|\le 1$, we have the inclusion
\begin{equation}\label{ee38}
  bC^{\a}(\R^{2})\subseteq bC^{\a}_{\mathcal{L}}(\R^{2}).
\end{equation}
Preferring simplicity to generality, we shall use H\"older spaces defined in terms of the
Euclidean norm (cf. Assumption \ref{Ass1}): by \eqref{ee38}, this results in a slightly more
restrictive condition compared to the analogous one given in terms of the homogeneous norm. On the
other hand, all the results of this paper can be proved using the homogeneous norm
$|x|^{\frac{1}{3}}+|v|$ as in \cite{Polidoro2}, \cite{DiFrancescoPascucci2} and \cite{Menozzi10}:
this would be more natural but would greatly increase the technicalities.

\smallskip We close this section by giving some standard Gaussian
estimates that play a central role in the parametrix construction. After the change of variables
\eqref{e39} with $\s=0$, the Langevin operator $\mathcal{L}_{a}$ is transformed into
\begin{equation}\label{e46}
 \LL_{a} =\frac{a}{2}\bar{\mathbf{V}}_{t}^{2}-\p_{t},\qquad \bar{\mathbf{V}}_{t}:=\p_{v}-t\p_{x}.
\end{equation}
Since $\LL_{a}$ is a heat operator with time dependent coefficients, its fundamental solution is
the Gaussian function $\pG_{a}(t,z;s,\z)=\pG_{a}(t,z-\z;s,0)$ where
\begin{equation}\label{ee40}
  \pG_{a}(t,x,y;s,0,0)=\frac{\sqrt{3}}{a\pi(t-s)^{2}}\exp\left(-\frac{2}{a(t-s)^{3}}\left(3x^{2}+3xy(t+s)+y^{2}(t^{2}+ts+s^{2})\right)\right)
\end{equation}
for $s<t$ and $x,y\in\R$.
\begin{lemma}\label{l3}
For every $\e>0$ there exists a positive constant $c$ such that
\begin{align}\label{ee41}
 \left|\bar{\mathbf{V}}_{t}\pG_{a}(t,x,y;s,0,0)\right|\le
 \frac{c}{\sqrt{t-s}}\pG_{a+\e}(t,x,y;s,0,0),\\ \label{ee42}
 \left|\bar{\mathbf{V}}^{2}_{t}\pG_{a}(t,x,y;s,0,0)\right|\le \frac{c}{t-s}\pG_{a+\e}(t,x,y;s,0,0),
\end{align}
for every $0\le s<t\le T$ and $x,y\in\R$.
\end{lemma}
\proof We remark that $\pG_{a}(t,x,y;s,0,0)$ has different asymptotic regimes as $t\to s^{+}$
depending on whether or not $s$ is zero: in fact, if $s=0$ then the quadratic form in the exponent
of $\pG_{a}$ is similar to that of the Langevin operator, that is
 $$\frac{1}{a}\begin{pmatrix}
   \frac{6}{t^{3}} & \frac{3}{t^{2}} \\
   \frac{3}{t^{2}} & \frac{2}{t} \
 \end{pmatrix}.$$
On the other hand, if $s\neq 0$ we see in \eqref{ee40} that all the components of the quadratic
form are $O((t-s)^{-3})$ as $t\to s^{+}$.

The thesis is a consequence of the following elementary inequality: for any $\e>0$ and $n\in\N$
there exists a positive constant $c_{\e,n}$ such that
\begin{equation}\label{e43}
  |\l|^{n} e^{-\frac{\l^{2}}{\m}}\le c_{n,\e}e^{-\frac{\l^{2}}{\m+\e}},\qquad \l\in\R.
\end{equation}
Indeed, we have
\begin{align}
  \left|\bar{\mathbf{V}}_{t}\pG_{a}(t,x,y;s,0,0)\right|&=\frac{1}{\sqrt{t-s}}\left|\frac{6 x+2v(t+2 s )}{a (t-s
  )^{\frac{3}{2}}}\right|\pG_{a}(t,x,y;s,0,0)\le
\intertext{(by \eqref{e43} with $n=1$)}
 &\le \frac{C}{\sqrt{t-s}}\pG_{a+\e}(t,x,y;s,0,0).
\end{align}
The proof of \eqref{ee42} is similar, using that
  $$\bar{\mathbf{V}}^{2}_{t}\pG_{a}(t,x,y;s,0,0)=\frac{4}{a(t-s)}\left(\frac{(3x+v(t+2s))^{2}}{a(t-s)^{3}}-1\right)\pG_{a}(t,x,y;s,0,0).$$
\endproof

\section{Assumptions and main results}\label{results}
We introduce the functional spaces used throughout the paper. For convenience, we give the
de\-finitions in the general multi-dimensional setting even if, except for Section \ref{pointw},
we will mainly consider dimension $d=2$.

Let $k,d\in\N$, $0<\a<1$ and $0\le t<T$. Denote by $m\B_{t,T}$ the space of all real-valued Borel
measurable functions $f=f_{s}(z)$ on $[t,T]\times \R^d$ and
\begin{itemize}
  \item $C^{\a}_{t,T}$ (resp. $bC^{\a}_{t,T}$) is the space of (resp. bounded) functions $f\in m\B_{t,T}$
  that are $\a$-H\"older continuous in $z$ uniformly with respect to $s$, that is
  $$\sup_{s\in[t,T]\atop z\neq \z}\frac{|f_{s}(z)-f_{s}(\z)|}{|z-\z|^{\a}}<\infty;$$
  \item $C^{k+\a}_{t,T}$ (resp. $bC^{k+\a}_{t,T}$) is the space of functions $f\in m\B_{t,T}$ that are $k$-times differentiable with respect to $z$
  with derivatives in $C^{\a}_{t,T}$  (resp. $bC^{\a}_{t,T}$).
\end{itemize}
We use boldface to denote the stochastic version of the previous
functional spaces. More precisely, let $\mathcal{P}_{t,T}$ be the predictable $\s$-algebra on
$[t,T]\times\O$.
\begin{definition}
We denote by $\mathbf{C}^{k+\a}_{t,T}$ the family of functions $f=f_{s}(z,\o)$ on
$[t,T]\times\R^{d}\times\O$ such that:
\begin{itemize}
 \item[i)] $(z,x)\mapsto f_{s}(z,\o)\in C^{k+\a}_{t,T}$ for any $\o\in\O$;
 \item[ii)] $(s,\w)\mapsto f_{s}(z,\o)$ is $\mathcal{P}_{t,T}$-measurable for any $z\in\R^{d}$.
\end{itemize}
Similarly, we define $\mathbf{bC}^{k+\a}_{t,T}$.
\end{definition}

\begin{definition}\label{d1} A stochastic fundamental solution $\mathbf{\pG}=\mathbf{\pG}(t,x,v;\t,\x,\y)$ for
the SPDE \eqref{spde} is a function defined for $0\le\t<t\le T$ and $x,v,\x,\y\in\R$, such that
for any $(\t,\z)\in[0,T)\times\R^{2}$ we have:
\begin{itemize}
  \item[i)] $\mathbf{\pG}(\cdot,\cdot,\cdot;\t,\z)$ belongs to $\mathbf{C}_{t_{0},T}(\R^{2})$,
  is twice continuously differentiable in $v$ and with probability one satisfies
\begin{equation}\label{e50}
\begin{split}
 \mathbf{\pG}(t,\gY_{t}(x,v);\t,\z)=\mathbf{\pG}(t_{0},x,v;\t,\z)
 &+\int_{t_{0}}^{t}\frac{1}{2}a_{s}(\gY_{s}(x,v))\left(\p_{vv}\mathbf{\pG}\right)(s,\gY_{s}(x,v);\t,\z)ds\\
 &+
 \int_{t_{0}}^{t}\s_{s}(\gY_{s}(x,v))\left(\p_v \mathbf{\pG}\right)(s,\gY_{s}(x,v);\t,\z)dW_{s}
\end{split}
\end{equation}
 for $\t<t_{0}\le t\le T$ and $x,v\in\R$, with $\gY=\gY_{t}(x,v)$ as in \eqref{e4};
  \item[ii)] for any bounded and continuous function $\phi$ on $\R^{2}$ and $z_{0}\in\R^{2}$, we have
    $$\lim_{(t,z)\to(\t,z_{0})\atop t>\t}\int_{\R^{2}}
    \mathbf{\pG}(t,z;\t,\z)\phi(\z)d\z=\phi(z_{0}),\qquad P\text{-a.s.}$$
\end{itemize}
\end{definition}

Next we state the standing assumptions on the coefficients of the SPDE \eqref{spde}.
\begin{assumption}[\bf Regularity]\label{Ass1}
$a\in\mathbf{bC}^{\a}_{0,T}$ for some $\a\in(0,1)$  and $\s\in\mathbf{bC}^{3+\a}_{0,T}$.
\end{assumption}
\begin{assumption}[\bf Coercivity]\label{Ass2}
There exists a random, finite and positive constant $\MM$ such that
\begin{equation}
 a_{t}(z)-\s_{t}^{2}(z)\ge \MM,\qquad t\in [0,T],\ z\in \R^{2},\ P\text{-a.s.}
\end{equation}
\end{assumption}
One of the main tools in our analysis is the following It\^o-Wentzell transform: for $\t\in[0,T)$
and $(x,v)\in\R^{2}$, we consider the one-dimensional SDE
\begin{equation}\label{sde}
 \gIW_{t,\t}(x,v)=v-\int_{\t}^t\s_s(x,\gIW_{s,\t}(x,v))dW_s.
\end{equation}
Assumption \ref{Ass1} ensures that \eqref{sde} is solvable in the strong sense and the map
$(x,v)\mapsto \left(x,\gIW_{t,\t}(x,v)\right)$ is a stochastic flow of diffeomorphisms of $\R^{2}$
(see Theorem \ref{t2} below). In Section \ref{ItoWen} we use this change of coordinates to
transform the SPDE \eqref{spde} into a PDE with random coefficients whose properties depend on the
gradient of the stochastic flow: to have a control on it, we impose the following additional
\begin{assumption} \label{Ass3}
There exist $\e>0$ and a random variable $M\in L^{p}(\O)$, with
$p>\max\left\{2,\frac{1}{\e}\right\}$, such that with probability one
\begin{align}
  \sup_{t\in[0,T]\atop (x,v)\in\R^2}(1+x^2+v^2)^{\e}|\p_{x}^{\b_{1}}\p_{v}^{\b_{2}}\s_t(x,v)|\leq M, \qquad &\b_{1}+\b_{2}=1, \\
  \sup_{t\in[0,T]\atop (x,v)\in\R^2}(1+x^2+v^2)^{\frac{1}{2}+\e}|\p_{x}^{\b_{1}}\p_{v}^{\b_{2}}
  \s_t(x,v)|\leq M, \qquad
  &\b_{1}+\b_{2}=2,3.
\end{align}
\end{assumption}
Assumption \ref{Ass3} is the main ingredient in the estimates of Section \ref{pointw}: it requires
that $\s_{t}(x,v)$ flattens as $(x,v)\to \infty$. In particular, this condition is clearly
satisfied if $\s=\s_{t}$ depends only on $t$ or, more generally, if the spatial gradient of $\s$
has compact support.

\smallskip We are now in position to state the main result of the paper.
\begin{theorem}\label{t1} Let Assumptions \ref{Ass1}, \ref{Ass2} and \ref{Ass3} be satisfied.
Then the Fokker-Planck SPDE \eqref{spde} has a fundamental solution $\mathbf{\pG}$ such that, for
some positive random variables $\m_{1}$ and $\m_{2}$, with probability one we have
%
\begin{align}\label{e31bbis}
\m^{-1}\pG^{\text{\rm heat}}\left( \m^{-1} Q_{t-\t}, g^{\text{\rm\tiny
IW},-1}(z)-\gvar_{t}^{\t,\z}\right) &\leq \mathbf{\pG}(t,z;\t,\z) \leq \m\pG^{\text{\rm
heat}}\left( \m Q_{t-\t}, g^{\text{\rm\tiny IW},-1}(z)-\gvar_{t}^{\t,\z}\right),\\ \label{e32bbis}
 \left| \p_{v}\mathbf{\pG}(t,x,v;\t,\z)\right|&\leq
 \frac{\m}{\sqrt{t-\t}}\pG^{\text{\rm heat}}\left( \m Q_{t-\t}, g^{\text{\rm\tiny IW},-1}(z)-\gvar_{t}^{\t,\z}\right),\\ \label{e33bbis}
 \left| \p_{vv}\mathbf{\pG}(t,x,v;\t,\z)\right|
 &\leq \frac{\m}{t-\t}\pG^{\text{\rm heat}}\left( \m Q_{t-\t}, g^{\text{\rm\tiny IW},-1}(z)-\gvar_{t}^{\t,\z}\right),
\end{align}
for every $0\leq\t<t\leq T$ and $z=(x,v),\z\in \R^2$, where
$g^{\text{\rm\tiny IW},-1}$ denotes the inverse of the stochastic flow
$(x,v)\mapsto \left(x,\gIW_{t,\t}(x,v)\right)$ defined by \eqref{sde} and
$\gvar_{t}^{\t,\z}$ is the integral curve (see Theorem \ref{t2} below), starting from $\z$, of
the vector field
\begin{equation}
  \Qg_{t,\t}=\left(\gIW_{t,\t},-\frac{\gIW_{t,\t} \p_x\gIW_{t,\t}}{\p_v\gIW_{t,\t}}\right),
\end{equation} $Q_{t}$ is defined as in \eqref{e7c} and
$$\G^{\text{\rm heat}}(A,z)=\frac{1}{2\pi\sqrt{\det A}}e^{-\frac{1}{2}\langle A^{-1}z,z \rangle}$$
is the two-dimensional Gaussian kernel with symmetric and positive definite covariance matrix $A$.
%
\end{theorem}
The proof of Theorem \ref{t1} 
is postponed to Section \ref{proof}.

\section{Pointwise estimates for It\^o processes}\label{pointw}
In this section we prove some estimate for stochastic flows of diffeomorphisms that will play a
central role in our analysis. Information about stochastic flows in a more general framework can
be found in \cite{MR1070361}. Since the following results are of a general nature and may be of
independent interest, in this section we reset the notations and give the proofs in the more
general multi-dimensional setting.

Specifically, until the end of
the section, the point of $\R^{d}$ is denoted by $z=(z_{1},\dots,z_{d})$ and we set
$\nabla_{z}=(\p_{z_{1}},\dots,\p_{z_{d}})$ and 
$\p^{\b}=\p_{z}^{\b}=\p_{z_{1}}^{\b_{1}}\cdots \p_{z_{d}}^{\b_{d}}$ for any multi-index $\b$.
We will also employ the
notation
 $$\langle z\rangle:=\sqrt{1+|z|^2},\qquad z\in\R^{d}.$$

First, we recall some basic facts about stochastic flows of diffeomorphisms. Let $k\in\N$. A
$\R^d$-valued random field $\phi_{\t,t}(z)$, with $0\le \t\le t\le T$ and $z\in\R^{d}$, defined on
$(\O,\F,P)$, is called a (forward) stochastic flow of $C^k$-diffeomorphisms if there exists a set
of probability one where:
\begin{itemize}
\item[i)] $\phi_{t,t}(z)=z$ for any $t\in [0,T]$ and $z\in\R^d$;
\item[ii)] $\phi_{\t,t}=\phi_{s,t}\circ\phi_{\t,s}$ for $0\le \t\le s\le t\le T$;
\item[iii)] $\phi_{\t,t}:\R^{d}\longrightarrow\R^{d}$ is a $C^k$-diffeomorphism for all $0\le \t\le t\le T$.
\end{itemize}
Stochastic flows can be constructed as solutions of stochastic differential equations. Let $B$ a
$n$-dimensional Brownian motion and consider the stochastic differential equation
\begin{equation}\label{sde0}
 \phi_{t}(z)=z+\int_{\t}^{t}b_{s}(\phi_{s}(z))ds+\int_{\t}^{t}\s_s(\phi_s(z))dB_s
\end{equation}
where $b=(b^i_t(z))$, $\s=(\s^{ij}_t(z))$ are a $d$-valued and $(d\times n)$-valued processes
respectively, on $[0,T]\times\R^d\times \O$. The following theorem summarizes the results of
Lemmas 4.5.3-7 and Theorems 4.6.4-5 in \cite{MR1070361}.
\begin{theorem}\label{t2} Let $b$, $\s\in \mathbf{bC}^{k,\a}_{0,T}$ for some $k\in\N$ and $\a\in(0,1)$. Then the
solution of the stochastic differential equation \eqref{sde0} has a modification $\phi_{\t,t}$
that is a forward stochastic flow of $C^k$-diffeomorphisms. Moreover,
$\phi_{\t,\cdot}\in\mathbf{C}_{\t,T}^{k,\a'}$ for any $\a'\in[0,\a)$ and $\t\in [0,T)$, and we
have the following estimates: for each $p\in\R$ there exists a positive constant
$\mathbf{c}_{1,p}$ such that
\begin{equation}
 E\left[\langle\phi_{\t,t}(z)\rangle^{p} \right]\leq \mathbf{c}_{1,p} \langle z\rangle^p, \qquad z\in\R^{d},\label{E4}\\
\end{equation}
and for each $p\ge 1$ there exists a positive constant $\mathbf{c}_{2,p}$ such that
\begin{equation}
 E\left[\left| \p^{\b}\phi_{\t,t}(z)\right|^{p} \right]\leq \mathbf{c}_{2,p}, \qquad z\in\R^{d},\ p\ge 1, \ 1\le|\b|\le k. \label{E5}
\end{equation}
\end{theorem}

\smallskip Now, consider $\phi_{\t,t}$ as in Theorem \ref{t2},
$F_{i}=F_{i,t}(z;\z)\in\mathbf{C}^k_{0,T}(\R^{2d})$, $i=1,2,$ and a real Brownian motion $W$. The
goal of this section is to prove some pointwise estimate for the It\^o process
\begin{equation}\label{e13}
 I_{\t,t}(z):=\int_{\t}^{t}F_{1,s}(z;\phi_{\t,s}(z))dW_s+\int_{\t}^{t}F_{2,s}(z;\phi_{\t,s}(z))ds,\qquad 0\le\t\le t\le T,\ z\in\R^{d},
\end{equation}
in terms of the usual H\"older norm in $\R^{d}$
 $$|f|_{\a}=\sup_{z\in\R^{d}}|f(z)|+\sup_{z,\z\in\R^{d}\atop z\neq \z}\frac{|f(z)-f(\z)|}{|z-\z|^{\a}},\qquad \a\in(0,1),$$
under the following
\begin{assumption} \label{Ass4}
There exist $\e_1,\e_2\in\R$ with $\e:=\e_1+\e_2>0$ and a random variable $M\in L^{\bar{p}}(\O)$,
with $\bar{p}>\left(2\vee d\vee \frac{d}{\e}\right)
$, such that
 $$\sum_{|\b|\le k}\sup_{t\in[0,T]\atop z,\z\in\R^{d}} \langle z\rangle^{\e_1}\langle \z\rangle^{\e_2}
 |\p^{\b}_{z,\z}F_{i,t}(z;\z)|\le M \qquad
 i=1,2,\ P\text{-a.s.}$$
\end{assumption}

The main result of this section is the following theorem which provides global-in-space pointwise
estimates for the process in \eqref{e13}.
\begin{theorem}\label{est1}
Let $\phi_{\t,t}$ be as in Theorem \ref{t2} and $F^{(i)}\in\mathbf{C}^k_{0,T}(\R^{2d})$, $i=1,2$,
for some $k\in\N$. Let $I=I_{\t,t}(z)$ be as in \eqref{e13} and set
  $$I_{\t,t}^{(\d)}(z):=\langle z\rangle^{\d}I_{\t,t}(z).$$
Under Assumption \ref{Ass4}, for any $p,\a$ and $\d$ such that
 $$\left(2\vee d\vee \frac{d}{\e}\right)<p<\bar{p},\qquad 0\le\a<\frac{1}{2}-\frac{1}{p},\qquad 0\le \d<\e-\frac{d}{p},$$
there exists a (random, finite) constant $\mathbf{m}$ such that
\begin{equation}\label{E6}
 \sum_{|\b|\le k-1}|\p^{\b}I^{(\d)}_{\t,t}|_{1-\frac{d}{p}}\le \mathbf{m}(t-\t)^{\a} \qquad P\text{-a.s.}
\end{equation}
\end{theorem}
\proof The proof is based on a combination of sharp $L^{p}$-estimates, Kolmogorov continuity
theorem in Banach spaces and Sobolev embedding theorem.

Let us first consider the case  $k=1$. We prove some preliminary $L^p$-estimates for $I_{\t,t}$
and $\p^{\b} I_{\t,t}$ with $|\b|=1$. Below we denote by $\bar{c}$ various positive constants that
depend only on $p,d,T$ and the flow $\phi$. By Burk\"older's inequality we have
\begin{align}
 E\left[|I^{(\d)}_{\t,t}(z)|^p\right]&\leq \bar{c}\,\langle z\rangle^{\d p} E\left[\left( \int_{\t}^t
 F_{1,s}^2(z;\phi_{\t,s}(z))ds\right)^{\frac{p}{2}}\right]+
 \bar{c}\,\langle z\rangle^{\d p} E\left[\left( \int_{\t}^t
 F_{2,s}(z;\phi_{\t,s}(z))ds\right)^{p}\right]
 \le \intertext{(by H\"older's inequality)}
 &\leq \bar{c}\,\langle z\rangle^{\d p} (t-\t)^{\frac{p-2}{2}}\int_{\t}^t
 E\left[|F_{1,s}(z;\phi_{\t,s}(z))|^p \right]ds\\
 &\quad+
 \bar{c}\,\langle z\rangle^{\d p} (t-\t)^{p-1}\int_{\t}^t E\left[|F_{2,s}(z;\phi_{\t,s}(z))|^p \right]ds \le
\intertext{(by Assumption \ref{Ass4})}
 &\leq \bar{c}\,\langle z\rangle^{(\d-\e_{1}) p} (t-\t)^{\frac{p-2}{2}} \int_{\t}^t E\left[M^p \langle
 \phi_{\t,s}(z)\rangle^{-\e_2 p}\right] ds\le
\intertext{(by H\"older's inequality with conjugate exponents $q:=\frac{\bar{p}}{p}$ and $r$)
}
 &\leq \bar{c}\, \langle z\rangle^{(\d-\e_{1})p}(t-\t)^{\frac{p-2}{2}}\|M\|_{L^{\bar{p}}(\O)}^p\int_{\t}^t E\left[
 \langle\phi_{\t,s}(z)\rangle^{-\e_2pr}\right]^{\frac{1}{r}}ds\le \intertext{(by \eqref{E4})}
 &=\bar{c}\,\langle z\rangle^{(\d-\e) p}(t-\t)^{\frac{p}{2}}.
\label{emb1}
\end{align}
The same estimate holds for the gradient of $I_{\t,t}$, that is
\begin{equation}\label{emb2}
 E\left[|\nabla I^{(\d)}_{\t,t}(z)|^p\right]\leq \bar{c}\,\langle z\rangle^{(\d-\e) p}(t-\t)^{\frac{p}{2}}.
\end{equation}
Indeed, let us consider for simplicity only the case $\d=0$ since the general case is a
straightforward consequence of the product rule: for $j=1,\dots,d$, we have
\begin{align}
 E\left[|\p_{z_{j}} I_{\t,t}(z)|^p\right]&\leq \bar{c}E\left[\left|\int_{\t}^t
 \Big((\p_{z_{j}}F_{1,s})(z;\phi_{\t,s}(z))+
 \langle\nabla_{\z}F_{1,s}(z;\phi_{\t,s}(z)),\p_{z_{j}}\phi_{\t,s}(z)\rangle\Big)
 dW_s\right|^p\right]\\
 &\quad + \bar{c}E\left[\left|\int_{\t}^t
 \Big((\p_{z_{j}}F_{2,s})(z;\phi_{\t,s}(z))+
 \langle\nabla_{\z}F_{2,s}(z;\phi_{\t,s}(z)),\p_{z_{j}}\phi_{\t,s}(z)\rangle\Big)
 ds\right|^p\right]\\
 &\leq \bar{c}(t-\t)^{\frac{p-2}{2}}\sum_{i=1}^2\int_{\t}^t
 E\left[|(\p_{z_{j}}F_{i,s})(z;\phi_{\t,s}(z))|^p+\left|\langle\nabla_{\z}F_{i,s}(z;\phi_{\t,s}(z)),\p_{z_{j}}
 \phi_{\t,s}(z)\rangle\right|^p\right]ds.
\end{align}
The terms containing $\p_{z_{j}}F_{i,s}$ can be estimated as before, by means of Assumption
\ref{Ass4}. On the other hand, by H\"older's inequality with conjugate exponents $q$ and $r$ with
$1<q<\frac{\bar{p}}{p}$, for every $i,j=1,\cdots,d $ we have
\begin{align*}
 E\left[\left|\langle\nabla_{\z}F_{i,s}(z;\phi_{\t,s}(z)),
 \p_{z_{j}}\phi_{\t,s}(z)\rangle\right|^p\right]&\leq
 E\left[\left|\nabla_{\z}F_{i,s}(z;\phi_{\t,s}(z))\right|^{pq}\right]^{\frac{1}{q}}
 E\left[\left|\p_{z_{j}}\phi_{\t,s}(z)\right|^{pr}\right]^{\frac{1}{r}}\le
\intertext{(by Assumption \ref{Ass4} and \eqref{E5})}
 & \leq \bar{c}_{2,pr}^{\frac{1}{r}}
 E\left[M^{pq}\langle\phi_{\t,s}(z)\rangle^{-\e_2 pq}\right]^{\frac{1}{q}}\langle z\rangle^{-\e_1 p}\le
\intertext{(by H\"older inequality with conjugate exponents $\bar{q}:=\frac{\bar{p}}{pq}>1$ and
$\bar{r}$)}
 & \leq \bar{c}_{2,pr}^{\frac{1}{r}}
 \|M\|_{L^{\bar{p}}(\O)}^p  E\left[\langle \phi_{\t,s}(z)\rangle^{-\e_2 pq\bar{r}}\right]^{\frac{1}{q\bar{r}}}
 \langle z\rangle^{-\e_1 p}\le
\intertext{(by \eqref{E4})}
 & \leq \bar{c} \|M\|_{L^{\bar{p}}(\O)}^p\, \langle z\rangle^{-\e p}.
\end{align*}
This proves \eqref{emb2} with $\d=0$.

Now, we have
\begin{align}
 E\left[\| I^{(\d)}_{\t,t} \|^p_{W^{1,p}(\R^d)}\right]&=E\left[ \int_{\R^d}\left(|I^{(\d)}_{\t,t}(z)|^p+
 |\nabla I^{(\d)}_{\t,t}(z)|^p\right)dz\right]\le
\intertext{(by \eqref{emb1} and \eqref{emb2})}
 &\le \bar{c} (t-\t)^{\frac{p}{2}}
 \int_{\R^d}\langle z\rangle^{(\d-\e) p}dz=
\intertext{(since $(\e-\d) p>d$ by assumption)}
 &=\bar{c} (t-\t)^{\frac{p}{2}}.\label{e14}
\end{align}
Estimate \eqref{e14} and Kolmogorov's continuity theorem for processes with values in the Banach
space $W^{1,p}(\R^d)$ (see, for instance, \cite{MR1070361}, Theor.1.4.1) yield
\begin{equation}\label{emb3}
 \| I^{(\d)}_{\t,t} \|_{W^{1,p}(\R^d)}\le \mathbf{m} (t-\t)^{\a}, \qquad 0\le\t\le t\le T,\ P\text{-a.s.}
\end{equation}
for some positive and finite random variable $\mathbf{m}$ and for $\a\in[0,\frac{p-2}{2p})$. This
is sufficient to prove \eqref{E6} with $k=1$: in fact, by the Sobolev embedding theorem, we have
the following estimate of the H\"older norm
\begin{equation}\label{E6b}
 |I^{(\d)}_{\t,t}|_{1-\frac{d}{p}}\le N\| I^{(\d)}_{\t,t} \|_{W^{1,p}(\R^d)}
\end{equation}
where $N$ is a positive constant that depends only on $p$ and $d$. Thus, combining \eqref{E6} and
\eqref{E6b}, we get the thesis with $k=1$.

Noting that
\begin{equation}\label{e12}
\begin{split}
 \p_{z_{j}} I_{\t,t}(z)&=\int_{\t}^t \Big((\p_{z_{j}}F_{1,s})(z;\phi_{\t,s}(z))+
 \langle\nabla_{\z}F_{1,s}(z;\phi_{\t,s}(z)),\p_{z_{j}}\phi_{\t,s}(z)\rangle\Big)
 dW_s\\
 &\quad + \int_{\t}^t
 \Big((\p_{z_{j}}F_{2,s})(z;\phi_{\t,s}(z))+
 \langle\nabla_{\z}F_{2,s}(z;\phi_{\t,s}(z)),\p_{z_{j}}\phi_{\t,s}(z)\rangle\Big) ds,
\end{split}
\end{equation}
for $j=1,\dots,d$, the thesis with $k=2$ can be proved repeating the previous arguments and using
\eqref{E6} for $k=1$ and Assumption \ref{Ass4} with $k=2$.

We omit the complete proof for brevity and since in the rest of the paper we will use \eqref{E6}
only for $k=1,2$. The general result can be proved by induction, using the multi-variate Fa\`a di
Bruno's formula.\endproof

\begin{remark}\label{r1}
Let $I_{\t,t}$ as in \eqref{e13} with coefficients $\tilde{F}_{1},\tilde{F}_{2}\in
b\mathbf{C}^{1}_{0,T}(\R^{2d})$ and let $\d>0$ and $\a\in[0,\frac{1}{2})$. Applying Theorem
\ref{est1} with $F_{i,t}(z;\z):=\langle z \rangle^{-\d}\tilde{F}_{i,t}(z;\z)$, $i=1,2$, we get the
existence of a (random, finite) constant $\mathbf{m}$ such that, with probability one,
 $$|I_{\t,t}(z)|\le \mathbf{m}\,\langle z \rangle^{\d}(t-\t)^{\a}, \qquad 0\le \t\le t\le T,\ z\in\R^{d}.$$
\end{remark}

\section{It\^o-Wentzell change of coordinates}\label{ItoWen}
We go back to the main SPDE \eqref{spde} and suppose that Assumptions \ref{Ass1}, \ref{Ass2} and
\ref{Ass3} are satisfied. In this section we study the properties of a random change of variables
which plays the same role as transformation \eqref{e15} in Step 1 of Section \ref{sub01} for the
Langevin SPDE. The main result of this section is Theorem \ref{th1} which shows that this change
of variables transforms SPDE \eqref{spde} into a PDE with random coefficients.

We denote by $(x,\gIW_{t,\t}(x,v))$ the stochastic flow of diffeomorphisms of $\R^{2}$ defined by
equation \eqref{sde}, that is
\begin{equation}\label{sde1}
  \gIW_{t,\t}(x,v)=v-\int_{\t}^t\s_s(x,\gIW_{s,\t}(x,v))dW_s,\qquad 0\le\t\le t\le T,\ (x,v)\in\R^{2}.
\end{equation}
By Theorem \ref{t2}, $\gIW_{\cdot,\t}\in \mathbf{C}_{\t,T}^{3,\a'}$ for any $\a'\in [0,\a)$.
Global estimates for $\gIW$ and its derivatives are provided in the next:
\begin{lemma}\label{lemma1} There exists
$\e\in\left(0,\frac{1}{2}\right)$ and a (random, finite) constant $\mathbf{m}$ such that, with
probability one,
\begin{align}\label{e20}
 |\gIW_{t,\t}(x,v)|&\le \mathbf{m} \sqrt{1+x^2+v^2},\\ \label{e21}
 e^{-\mathbf{m}(t-\t)^{\e}}\le \p_v\gIW_{t,\t}(x,v)&\le e^{\mathbf{m}(t-\t)^{\e}},\\ \label{e22}
 |\p_x \gIW_{t,\t}(x,v)|&\le \mathbf{m}(t-\t)^{\e},\\ \label{e23}
 |\p^{\b}\gIW_{t,\t}(x,v)|&\le \frac{\mathbf{m}(t-\t)^{\e}}{\sqrt{1+x^2+v^2}},
\end{align}
for any $(x,v)\in\R^{2}$, $0\le\t\le t\le T$ and $|\b|=2$.

\end{lemma}
\proof Estimate \eqref{e20} follows directly from Remark \ref{r1} (with $\d=1$). Differentiating
\eqref{sde1}, we find that $\p_{v}\gIW_{t,\t}$ solves the linear SDE
\begin{equation}
  \p_v\gIW_{t,\t}(x,v)=1-\int_{\t}^{t}(\p_2\s_s)(x,\gIW_{s,\t}(x,v))\p_v\gIW_{s,\t}(x,v)dW_s,
\end{equation}
where $\p_2\s_t$ denotes the partial derivative of $\s_{t}(\cdot,\cdot)$ with respect to its
second argument. Hence we have
\begin{equation}\label{e01}
 \p_v\gIW_{t,\t}(x,v)=\exp\left(-\int_{\t}^{t}(\p_2\s_s)(x,\gIW_{s,\t}(x,v))dW_s-
 \frac{1}{2}\int_{\t}^{t}(\p_2\s_s)^2(x,\gIW_{s,\t}(x,v))ds\right).
\end{equation}
Now we apply Theorem \ref{est1} with $\phi_{\t,t}(x,v)=(x,\gIW_{t,\t}(x,v))$ and
$F_{i,t}(z;x,V)=(\p_2\s_t(x,V))^i$, $i=1,2$: thanks to Assumption \ref{Ass3}, we get estimate
\eqref{e21}. Incidentally, from Theorem \ref{est1} we also deduce that the first order derivatives
of $\p_v\gIW_t$ are bounded:
\begin{equation}\label{e21b}
  |\p^{\b}\p_v \gIW_{t,\t}(x,v)|\le \mathbf{m}(t-\t)^{\e},\qquad |\b|=1.
\end{equation}
This last estimate is used in the next step, for the proof of \eqref{e22}.

Similarly, we have
\begin{equation}
 \p_{x}\gIW_{t,\t}(x,v)=-\int_{\t}^{t}\big((\p_1\s_s)(x,\gIW_{s,\t}(x,v))+(\p_2\s_s)(x,\gIW_{s,\t}(x,v))\p_{x}\gIW_{s,\t}(x,v)\big)dW_s.
\end{equation}
Thus, we have a linear SDE whose solution is given by
\begin{align}
 \p_{x}\gIW_{t,\t}(x,v)&=-\p_{v}\gIW_{t,\t}(x,v)\int_{\t}^{t}\frac{(\p_1\s_s)(x,\gIW_{s,\t}(x,v))}{\p_v\gIW_{s,\t}(x,v)}dW_s\\
 &\quad-\p_{v}\gIW_{t,\t}(x,v)\int_{\t}^{t}\frac{(\p_1\s_s)(x,\gIW_{s,\t}(x,v))(\p_2\s_s)(x,\gIW_{s,\t}(x,v))}{\p_v\gIW_{s,\t}(x,v)}ds,
\end{align}
Again, \eqref{e22} follows from Theorem \ref{est1} thanks to Assumption \ref{Ass3} and estimates
\eqref{e21} and \eqref{e21b}.

Eventually, the same argument can be used to prove \eqref{e23}: indeed, differentiating
\eqref{sde1} we have that $\p^{\b}\gIW_t$ satisfies a linear SDE whose solution is explicit. 
Thus, for $|\b|=2$, $\p^{\b}\gIW_t$ can be expressed in the form \eqref{e13} with the coefficients
satisfying Assumption \ref{Ass4} for some $\e>1$. Applying Theorem \ref{est1} with $\d=1$ we get
estimate \eqref{e23}.
\endproof

Next, we provide a version of the It\^o-Wentzell formula for an equation of the form
\begin{equation}\label{eIW}
 d_{\mathbf{Y}}u_{t,\t}(x,v)=f_t(x,v)dt+g_t(x,v)dW_t,\qquad \mathbf{Y}=\p_{t}+v\p_{x},
\end{equation}
with $u,f,g \in \mathbf{C}_{\t,T}$. Equation \eqref{eIW} is solved in the strong sense which means
 $$u_{t,\t}\left(\gY_{t-\t}(x,v)\right)=u_{\t,\t}(x,v)+
 \int_{\t}^{t}f_{s}(\gY_{s-\t}(x,v))ds+\int_{\t}^{t}g_{s}(\gY_{s-\t}(x,v))dW_{s},\qquad t\in[\t,T],$$
with probability one, where $\gY_{t}(x,v)=(x+tv,v)$ is the integral curve in $\R^{2}$ of the
vector field $v\p_{x}$, starting from $(x,v)$. The following lemma shows how \eqref{eIW} is
modified by the It\^o-Wentzell transform
\begin{equation}\label{change}
  \hat{u}_{t,\t}(x,v)=u_{t,\t}(x,\gIW_{t,\t}(x,v)),
\end{equation}
with $\gIW_{t,\t}$ as in \eqref{sde1}.
\begin{lemma}[\bf It\^o-Wentzell formula]\label{lemma2}
Let $\p_{2}u_{t,\t}, \p_{22}u_{t,\t},\p_{2}g_{t} \in \mathbf{C}_{\t,T}$ and assume that
\eqref{eIW} holds. Then $\hat{u}_{t,\t}$ in \eqref{change} solves
\begin{equation}\label{eIW1}
 d_{\mathbf{\hat{Y}}}\hat{u}_{t,\t}(x,v)=F_{t}(x,v)dt+G_{t}(x,v)dW_t,
\end{equation}
where
\begin{align}\label{e24}
 F_{t}(x,v)&=\hat{f}_t(x,v)+\frac{1}{2}\hat{\s}^2_t(x,v)\widehat{\p_{22}u_{t,\t}}(x,v)-
 \widehat{\p_2g_t}(x,v)\hat{\s}_t(x,v),\\
 G_{t}(x,v)&=\hat{g}_t(x,v)-\hat{\s}_t(x,v)\widehat{\p_2u_{t,\t}}(x,v),\\ \label{e25}
 \mathbf{\hat{Y}}&=\p_t+\gIW_{t,\t}\p_x-\frac{\gIW_{t,\t} \p_x\gIW_{t,\t}}{\p_v\gIW_{t,\t}}\p_v.
\end{align}
Moreover, we have
\begin{align}\label{eIW1bis}
 \widehat{\p_2u_{t,\t}}=\frac{\p_v\hat{u}_{t,\t}}{\p_v\gIW_{t,\t}},\qquad
 \widehat{\p_{22}u_{t,\t}}=\frac{\p_{vv}\hat{u}_{t,\t}}{(\p_v\gIW_{t,\t})^{2}}-\frac{\p_{vv}\gIW_{t,\t}
 \p_v\hat{u}_{t,\t}}{(\p_v\gIW_{t,\t})^{3}}.
\end{align}
\end{lemma}
\proof We have to show that
\begin{align}
 \hat{u}_{t,\t}\left(\gvar_{t,\t}(x,v)\right)=\hat{u}_{0}(x,v)+
 \int_{0}^{t}F_{s}(\gvar_{s,\t}(x,v))ds+\int_{0}^{t}G_{s}(\gvar_{s,\t}(x,v))dW_{s},\qquad
 t\in[0,T],
\end{align}
where $\gvar_{t,\t}(x,v)$ is the integral curve, starting from $(x,v)$, of the vector field
\begin{equation}\label{e26}
  \Qg_{t,\t}=\left(\gIW_{t,\t},-\frac{\gIW_{t,\t} \p_x\gIW_{t,\t}}{\p_v\gIW_{t,\t}}\right).
\end{equation}
Notice that, with the usual identification of vector fields with first order operators, we have
$\mathbf{\hat{Y}}=\p_{t}+\Qg_{t,\t}$. Moreover, $\gvar$ is well defined thanks to the estimates of
Lemma \ref{lemma1}.

If $u\in \mathbf{C}^{2}_{\t,T}$ then \eqref{eIW} can be written in the usual It\^o sense
\begin{equation}\label{eIW2}
 du_{t,\t}(x,v)=(f_t(x,v)-v\p_xu_{t,\t}(x,v))dt+g_t(x,v)dW_t.
\end{equation}
Then, by the standard It\^o-Wentzell formula (see, for instance, Theor. 3.3.1 in
\cite{MR1070361}), we have
\begin{align}\label{eIW3}
 d\hat{u}_{t,\t}=\left((\hat{f}_t-\gIW_{t,\t}\widehat{\p_1u_{t,\t}})+\frac{1}{2}\hat{\s}^2_t\widehat{\p^2_2u_{t,\t}}-
 \widehat{\p_2g_t}\hat{\s}_t\right)dt+\left(\hat{g}_t-\hat{\s}_t\widehat{\p_2u_{t,\t}}\right)dW_t.
\end{align}
From the chain rule we easily derive equations \eqref{eIW1bis} and also
  $$\widehat{\p_1u_{t,\t}}=\p_x\hat{u}_{t,\t}-\frac{\p_x\gIW_{t,\t}}{\p_v\gIW_{t,\t}}\p_v\hat{u}_{t,\t}.$$
Plugging these formulas into \eqref{eIW3} we get \eqref{eIW1}.

In the general case, it suffices to
apply what we have just proved to a smooth approximation in $(x,v)$ of $u_{t,\t}$ and then pass to the limit.
\endproof

Applying the It\^o-Wentzell formula to SPDE \eqref{spde} we get the following
\begin{theorem}\label{th1}
Let $u_{t,\t},\p_{2}u_{t,\t},\p_{22}u_{t,\t}\in \mathbf{C}_{\t,T}$ and let Assumptions \ref{Ass1},
\ref{Ass2} and \ref{Ass3} be satisfied. If $u_{t,\t}$ solves the SPDE \eqref{spde} then
$\hat{u}_{t,\t}$ in \eqref{change} is such that
$\hat{u}_{t,\t},\p_{v}\hat{u}_{t,\t},\p_{vv}\hat{u}_{t,\t}\in \mathbf{C}_{\t,T}$ and
\begin{equation}\label{PDE}
 d_{\mathbf{\hat{Y}}}\hat{u}_{t,\t}(x,v)=\big(\au_{t,\t}(x,v)\p_{vv}\hat{u}_{t,\t}(x,v)+\bu_{t,\t}(x,v)\p_v\hat{u}_{t,\t}(x,v)\big)dt,
\end{equation}
with $\mathbf{\hat{Y}}$ as in \eqref{e25} and
\begin{align}\label{e04}
 \au_{t,\t}=\frac{\hat{a}_t-\hat{\s}^2_t}{2(\p_v\gIW_{t,\t})^{2}},\qquad
 \bu_{t,\t}=-\frac{1}{(\p_v\gIW_{t,\t})^{2}}\left(\hat{\s}_t\p_{v}\hat{\s_{t}}+\frac{\left(\hat{a}_t-\hat{\s}^2_t\right)
 \p_{vv}\gIW_{t,\t}}{2\p_v\gIW_{t,\t}}\right).
\end{align}
\end{theorem}

\proof The thesis follows from the It\^o-Wentzell formula of Lemma \ref{lemma2} with
$f_t=\frac{1}{2}a_t\p_{22}u_{t,\t}$ and $g_t=\s_t\p_2u_{t,\t}$: the assumptions $\p_{2}u_{t,\t},
\p_{22}u_{t,\t},\p_{2}g_{t} \in \mathbf{C}_{\t,T}$ are clearly satisfied.
\endproof

\section{Time-dependent parametrix method}\label{mproof}
In this section we study equation \eqref{PDE} for fixed $\o\in\O$ and $0\le\t<T<\infty$. More
generally, we consider a deterministic equation of the form
\begin{equation}\label{pde0}
 \mathcal{K}_tu_t(z)=\mathcal{L}_tu_t(z)-\p_tu_t(z)=0
\end{equation}
where 
\begin{equation}
 \mathcal{L}_tu_t(z):=\frac{1}{2}a_t(z)\p_{vv}u_t(z){+b_t(z)\p_{v}u_t(z)}-\langle Y_t(z),\nabla_z u_t(z) \rangle, \qquad t\in[\t,T],\ z=(x,v)\in\R^2,
\end{equation}
and $Y_{t}=(Y_{1,t},Y_{2,t})$ is a generic vector field. We assume the following conditions on the
coefficients.
\begin{assumption}\label{As1}
There exist positive constants $\a,\l_1$ such that $a, {b} \in C_{\t,T}^{\a}$ with H\"older
constant $\l_1$ and
\begin{equation}\label{as1}
 \l_1^{-1}\le a_t(z)\le \l_1, \quad {|b_t(z)|\le \l_1} \qquad (t,z)\in[\t,T]\times \R^2.
\end{equation}
\end{assumption}
\begin{assumption}\label{As2} $Y\in C_{\t,T}$ and is uniformly Lipschitz continuous in the sense that
\begin{equation}\label{lipc}
  \sup_{t\in[\t,T]\atop z\neq \z}\frac{|Y_{t}(z)-Y_{t}(\z)|}{|z-\z|}\le \l_{2}
\end{equation}
for some positive constant $\l_2$. Moreover $\p_vY_{1,t}\in C^{\a}_{\t,T}$ and
\begin{equation}\label{as2}
 \l_2^{-1}\le \p_vY_{1,t}(z)\le \l_2, \qquad (t,z)\in[\t,T]\times \R^2.
\end{equation}
\end{assumption}
\begin{remark} \label{r2} When the coefficients are smooth,
conditions \eqref{as1} and \eqref{as2} ensure the validity of the weak H\"ormander condition:
indeed the vector fields $\sqrt{a}\p_v$ and $Y$, together with their commutator, span $\R^{2}$ at
any point. In this case a smooth fundamental solution exists by H\"ormander's theorem.
\end{remark}

Since the coefficients are assumed to be only measurable in time, a solution to \eqref{pde0} has
to be understood in the integral sense according to the following definition. 
\begin{definition}\label{d1bis} A fundamental solution $\pG=\pG(t,z;s,\z)$ for
equation \eqref{spde} is a function defined for $\t\le s<t\le T$ and $z,\z\in\R^2$, such that for
any $(s,\z)\in[\t,T)\times\R^{2}$ we have:
\begin{itemize}
  \item[i)] for $s<t_{0}\le t\le T$ and $z\in\R^2$, $\pG(\cdot,\cdot;s,\z)$ belongs to $C_{t_{0},T}$,
  is twice continuously differentiable in $v$ and satisfies
\begin{equation}
 \mathbf{\pG}(t,\gvar_{t}^{t_0,z};s,\z)=\pG(t_{0},z;s,\z)
 +\int_{t_{0}}^{t}\left(\frac{1}{2}a_{\r}(\gvar_{\r}^{t_0,z})\left(\p_{vv}\pG\right)(\r,\gvar_{\r}^{t_0,z};s,\z)
 {+b_{\r}(\gvar_{\r}^{t_0,z})\left(\p_{v}\pG\right)(\r,\gvar_{\r}^{t_0,z};s,\z)}\right)d\r
\end{equation}
 where $\gvar_{t}^{t_0,z}$ stands for the integral curve of the field $Y$ with
 initial datum $\gvar_{t_0}^{t_0,z}=z$;
  \item[ii)] for any bounded and continuous function $\phi$ and $z_{0}\in\R^{2}$, we have
    $$\lim_{(t,z)\to(s,z_{0})\atop t>s}\int_{\R^{2}}
    \pG(t,z;s,\z)\phi(\z)d\z=\phi(z_{0}).$$
\end{itemize}
\end{definition}
The main result of this section is the following
\begin{theorem}\label{t3}
Let Assumptions \ref{As1} and \ref{As2} be in force. Then the PDE \eqref{pde0} has a fundamental
solution $\G$ such that, for any $z=(x,v),\z \in \R^2$ and $\t\le s< t\le T$,
\begin{align}\label{e40}
 \m^{-1}\G^{\text{\rm heat}}\left(\m^{-1}Q_{t-s},z-\gvar_{t}^{s,\z}\right)&\le \G(t,z;s,\z) \le
 \m\G^{\text{\rm heat}}\left(\m Q_{t-s},z-\gvar_{t}^{s,\z}\right),\\ \label{e41}
 \left|\p_v\G(t,x,v;s,\z)\right| &\le \frac{\m}{\sqrt{t-s}}\G^{\text{\rm heat}}\left(\m Q_{t-s},z-\gvar_{t}^{s,\z}\right),\\ \label{e42}
 \left|\p_{vv}\G(t,x,v;s,\z)\right| &\le \frac{\m}{t-s}\G^{\text{\rm heat}}\left(\m Q_{t-s},z-\gvar_{t}^{s,\z}\right).
\end{align}
where $Q_{t}$ is as in \eqref{e7c}, $\gvar_{t}^{s,\z}$ is as in Definition \ref{d1bis} and $\m$ is
a positive constant that depends only on $\l_1, \l_2$, {$\a$} and $T$.
\end{theorem}

\subsection{Proof of Theorem \ref{t3}}
We prove the Gaussian bounds \eqref{e40}-\eqref{e42} under the assumption that the coefficients
$a$ and $Y$ are smooth; by Remark \ref{r2}, this guarantees the existence of a smooth fundamental
solution. Our estimates extend and sharpen classical Gaussian bounds for H\"ormander operators
(e.g. \cite{MR865430}). Moreover, our estimates are independent of the regularity of the
coefficients and, by standard regularization arguments, lead to a priori Gaussian bounds for
operators satisfying Assumptions \ref{as1}-\ref{as2}.

\subsubsection{Parametrix expansion}
For fixed $(s,\y)\in[\t,T)\times\R^{2}$, let
\begin{equation}\label{e31}
 \gvar_{t}^{s,\y}=\y+\int_s^t Y_{\r}\left(\gvar_{\r}^{s,\y}\right)d\r,\qquad t\in[\t,T],
\end{equation}
be the integral curve of $Y_{t}$ starting from $(s,\y)$. Following \cite{MR2659772} we linearize
$Y_{t}=Y_{t}(z)$ at $(s,\y)$ setting
  $$\bar{Y}_{t}^{s,\y}(z)=
  Y_{t}(\gvar_{t}^{s,\y})+\left(D Y_{t}\right)(\gvar_{t}^{s,\y})
  \left(z-\gvar_{t}^{s,\y}\right),\qquad t\in[s,T],\ z\in\R^{2}.$$
where $DY_t$ stands for a reduced Jacobian defined as
 $$DY_t:=\begin{pmatrix}0 & \p_vY_{1,t} \\ 0 & 0
\end{pmatrix}.$$
Then we consider the linear approximation of $\mathcal{L}_{t}$ defined as
\begin{equation}\label{e30bis}
 \bar{\mathcal{L}}_{t}^{s,\y}:=\frac{1}{2}a_{t}(\gvar_{t}^{s,\y})\p_{vv}-
 \langle \bar{Y}_{t}^{s,\y}(z),\nabla\rangle.
\end{equation}
The diffusion coefficient of $\bar{\mathcal{L}}_{t}^{s,\y}$ depends on $t$ only (apart from $s,\y$
that are fixed parameters), while the drift coefficients depend on $t$ and linearly on $x,v$.
Notice that 
$\bar{\mathcal{L}}_{t}^{s,\y}-\p_t$ is the forward Kolmogorov operator of the system of linear
SDEs
\begin{equation}\label{e30tris}
 dH_{t}= \bar{Y}_{t}^{s,\y}\left(H_{t}\right)dt+
 \sqrt{a_{t}(\gvar_{t}^{s,\y})}\mathbf{e_2}dW_{t}.
\end{equation}
Let $H_t^{t_{0},\z}$ denote the solution of \eqref{e30tris} starting from $\z$ at time
$t_{0}\in[s,T)$. Then $H_t^{t_{0},\z}$ is a Gaussian process: the mean
$\bar{\g}^{s,\y}_{t,t_{0}}(\z):=E\left[H_t^{t_{0},\z}\right]$ solves the ODE
\begin{equation}\label{e31bis}
 \bar{\g}^{s,\y}_{t,t_{0}}(\z)=\z+\int_{t_{0}}^t\bar{Y}_{\r}^{s,\y}\left(\bar{\g}^{s,\y}_{\r,t_{0}}(\z)\right)
 d\r,\qquad t\in[t_{0},T],
\end{equation}
and the covariance matrix is given by
\begin{equation}\label{ee32}
 \AA_{t,t_{0}}^{s,\y}=\int_{t_{0}}^t
 a_{\r}(\gvar_{s}^{\r,\y})\left(E^{s,\y}_{t,\r}\mathbf{e}_2\right)
 \left(E^{s,\y}_{t,\r}\mathbf{e}_2\right)^{\ast}d\r,
 \end{equation}
where $E^{s,\y}_{t,\r}$ is the fundamental matrix associated with $(D
Y_{t})(\gvar_{t}^{s,\y})$, that is the solution of
\begin{equation}
  E^{s,\y}_{t,\r}=\text{Id}+\int_{\r}^{t}(D Y_{u})(\gvar_{u}^{s,\y})E^{s,\y}_{u,\r}du,\qquad
  t\in[\r,T],
\end{equation}
with $\text{Id}$ equal to the $(2\times 2)$-identity matrix.

\begin{lemma}\label{lemma3}
For any $\y\in\R^2$ and $\t\le s\le t_{0}<t\le T$, we have $\det\AA_{t,t_{0}}^{s,\y}>0$.
\end{lemma}
\proof
By Assumption \ref{As1} it is enough to prove the assertion for $a\equiv 1$.
Suppose that there exist $\z\in\R^2\setminus \left\{0\right\}$, $\y\in\R^2$ and $\t\le s\le t_0<t\le T$ such that
$\langle\AA_{t,t_0}^{s,\y}\z,\z\rangle = 0$.
Since $\AA_{t,t_0}^{s,\y}$ is positive semi-definite, this is equivalent to the condition
$$|(E^{s,\y}_{t,\r}\mathbf{e}_2)^{\ast}\z|^2=0, \quad \text{a.e.  } \r\in (t_0,t),$$
that is $((E^{s,\y}_{t,\r})^{\ast}\z)_2=0$, for a.e. $\r\in (t_0,t)$. We have
$$\p_{\r}(E^{s,\y}_{t,\r})^{\ast}\z=-D Y_{\r}^{\ast}(\gvar_{\r}^{s,\y})(E^{s,\y}_{t,\r})^{\ast}\z, $$
and therefore
\begin{equation}
0=\p_{\r}((E^{s,\y}_{t,\r})^{\ast}\z)_2=\p_v Y_{1,\r}(\gvar_{\r}^{s,\y})((E^{s,\y}_{t,\r})^{\ast}\z)_1.
\end{equation}
Since $((E^{s,\y}_{t,\r})^{\ast}\z)_2=0$ and $\p_v Y_{1,\r}\in [\l_2^{-1},\l_2]$ by Assumption
\ref{As2} we have $(E^{s,\y}_{t,\r})^{\ast}\z\equiv 0$,  for a.e. $\r\in (t_0,t)$, which is
absurd.
\endproof

Lemma \ref{lemma3} ensures that the Gaussian process in \eqref{e30tris} admits a transition
density that is the fundamental solution of $\bar{\mathcal{L}}_{t}^{s,\y}-\p_t$. To be more
precise, let us recall the notation $\G^{\text{\rm heat}}(A,z)$ for the two-dimensional Gaussian
kernel with covariance matrix $A$ (cf. Theorem \ref{t1}).
\begin{proposition}
For any $0\le\t\le s\le t_{0}<t\le T$ and $z,\z,\y\in\R^{2}$, the function
\begin{align}\label{fsol}
 \pG_{s,\y}(t,z;t_{0},\z):=\G^{\text{\rm heat}}
 \left(\AA_{t,t_{0}}^{s,\y},z-\bar{\g}^{s,\y}_{t,t_{0}}(\z)\right)
\end{align}
is the fundamental solution of $\bar{\mathcal{L}}_{t}^{s,\y}-\p_t$, evaluated at $(t,z)$ and with
pole at $(t_{0},\z)$.
\end{proposition}
We are now in position to define the parametrix $Z$ for $\mathcal{K}_{t}$ in \eqref{pde0}. We set
\begin{equation}
 Z(t,z;s,\z):=\pG_{s,\z}(t,z;s,\z),\qquad \t\le s<t\le T,\ z,\z\in\R^{2}.
\end{equation}
Since
   $$\g_{t}^{s,\z}=\z+\int_{s}^{t}Y_{\r}(\g_{\r}^{s,\z})d\r=\z+\int_{s}^{t}\bar{Y}_{\r}^{s,\z}(\g_{\r}^{s,\z})d\r$$
we have $\g_{t}^{s,\z}=\bar{\g}^{s,\z}_{t,s}(\z)$ and therefore the parametrix reads
\begin{equation}\label{parametrix}
 Z(t,z;s,\z)=\G^{\text{\rm heat}}\left(\AA_{t,s}^{s,\z},z-\gvar_{t}^{s,\z}\right)
\end{equation}
for $\t\le s<t\le T$ and $z,\z\in\R^2$.
The parametrix 
is an approximation of the fundamental solution $\pG$ of $\mathcal{K}_{t}$: 
indeed, since $Z(s,\cdot;s,\z)=\d_{\z}$ and $\G(t,z;t,\cdot)=\d_{z}$, we have
\begin{align}\label{expansion0}
 \pG(t,z;s,\z)-Z(t,z;s,\z)&=\int_{\R^2}\left(\pG(t,z;s,\y)Z(s,\y;s,\z)-\G(t,z;t,\y)Z(t,\y;s,\z)\right)d\y\\
 &=\int_{s}^t\int_{\R^2}-\p_{\r}\left(\pG(t,z;\r,\y)Z(\r,\y;s,\z)\right)d\y d\r\\
 &=\int_{s}^t\int_{\R^2}\left(\mathcal{L}^{\ast}_{\r}\pG(t,z;\r,\y)Z(\r,\y;s,\z)
 -\pG(t,z;\r,\y)\bar{\mathcal{L}}_{\r}^{s,\z}Z(\r,\y;s,\z)\right)d\y d\r=\\
 &=\int_{s}^t\int_{\R^2}\pG(t,z;\r,\y)(\mathcal{L}_{\r}-\bar{\mathcal{L}}_{\r}^{s,\z})Z(\r,\y;s,\z)d\y
 d\r=\\ &=\int_{s}^t\int_{\R^2}\pG(t,z;\r,\y)\mathcal{K}_{\r}Z(\r,\y;s,\z)d\y d\r.
\end{align}
Iterating the formula, for $N\ge 1$ we get the expansion
\begin{equation}\label{expansion}
\begin{split} \pG(t,z;s,\z)=Z(t,z;s,\z)&+\sum_{k=1}^{N-1}\int_{s}^t\int_{\R^2}Z(t,z;\r,\y)(\mathcal{K}_{\r}Z)_k(\r,\y;s,\z)d\y d\r\\
 &+\int_{s}^t\int_{\R^2}\G(t,z;\r,\y)(\mathcal{K}_{\r}Z)_N(\r,\y;s,\z)d\y d\r
\end{split}
\end{equation}
where
\begin{equation}\label{expansion2}
\begin{split}
 (\mathcal{K}_{t}Z)_1(t,z;s,\z)&=\mathcal{K}_{t}Z(t,z;s,\z)\\
 (\mathcal{K}_{t}Z)_{k+1}(t,z;s,\z)&=\int_{s}^t\int_{\R^2}\mathcal{K}_{t}Z(t,z;\r,\y)(\mathcal{K}_{\r}Z)_k(\r,\y;s,\z)d\y d\r
\end{split}
\end{equation}
As $N$ tends to infinity we formally obtain a representation of $\pG$ as a series of convolution
kernels. Unfortunately, as already noticed in \cite{MR2659772}, such an argument cannot be made
rigorous because of the transport term. The problem is that, using only the Gaussian estimates for
the parametrix, it seems difficult to control the iterated kernels uniformly in $k$.

For this reason, we first prove some bound for expansion \eqref{expansion} and estimate the
remainder via stochastic control techniques as in \cite{MR2659772}. Once we have obtained the
Gaussian bounds for the fundamental solution $\pG$, a posteriori we prove the convergence of the
series and the bounds for the derivatives of $\pG$.

\subsubsection{Gaussian bounds for the parametrix}\label{pbounds}
\begin{proposition}\label{p0}
There exists a positive constant $c$, only dependent on 
$\l_1$, $\l_2$ and $T$, such that
\begin{equation}\label{e36}
 c^{-1}|\mathcal{D}_{\sqrt{t-s}}z|^2\le \langle\AA_{t,s}^{s,\z}z,z \rangle\le c|\mathcal{D}_{\sqrt{t-s}}z|^2, \qquad \t\le s<t\le T,\
 z,\z\in\R^2,
\end{equation}
where, for $\l>0$, $\mathcal{D}_{\l}$ is the diagonal matrix $\text{diag}(\l^{3},\l)$ that is
the spatial part of the ultra-parabolic dilation operator \eqref{egdil}.
\end{proposition}
\proof
By Assumptions \ref{As1}  it is enough to prove the assertion for $a\equiv 1$.
For $\l>0$, let $\mathcal{U}_{\l}$ be the set of $2\times 2$, time-dependent matrices
$\mathcal{Y}_{t}$, with entries uniformly bounded by $\l$, and such that $(\mathcal{Y}_t)_{1,2}\in
[\l^{-1},\l]$. Let $\mathcal{Y}_t\in\mathcal{U}_{\l}$ and
 $$\mathcal{A}_{t,s}:=\int_s^t\left(\mathcal{E}_{t,\r}\mathbf{e_2}\right)
 \left(\mathcal{E}_{t,\r}\mathbf{e_2}\right)^{\ast}d\r, \qquad \t\le s< t\le T,$$
where $\mathcal{E}_{t,\r}$ denotes the resolvent associated with $\mathcal{Y}_t$. We split the
proof in two steps.

\textit{Step 1.} First we prove that
\begin{equation}\label{e34}
 c^{-1}|z|^2\le \langle\mathcal{A}_{1,0}z,z \rangle\le c|z|^2,
\end{equation}
where $c$ is a positive constant which depends only on $\l$. 
As in \cite{MR2659772} (see Proposition 3.4), we consider the map
 $$\Psi: L^2([0,1],\mathcal{M}_2(\R))\longrightarrow \R, \qquad \Psi(\mathcal{Y}):=\det \mathcal{A}_{1,0},$$
where $\mathcal{M}_2(\R)$ is the space of $2\times2$ matrices with real entries. Notice that
$\mathcal{U}_{\l}$ is compact in the weak topology of $L^2([0,1],\mathcal{M}_2(\R))$ because it is
bounded, convex and closed in the strong topology (cf., for instance, \cite{Brezis}, Corollary
III.19). On the other hand, $\Psi$ is continuous from $L^2([0,1],\mathcal{M}_2(\R))$, equipped
with the weak topology, to $\R$. Therefore the image $\Psi(\mathcal{U}_{\l})$ is a compact subset
of $\R_{>0}$ by Lemma \ref{lemma3}. Thus there exists $\bar{\l}>0$ such that
$\inf\{\det\mathcal{A}_{1,0}\mid \mathcal{Y}\in \mathcal{U}_{\l}\}\ge {\bar{\l}^{-1}}$ and
$\sup\{\|\mathcal{A}_{1,0}\|\mid \mathcal{Y}\in \mathcal{U}_{\l}\}\le {\bar{\l}}$. This suffices
to prove \eqref{e34}.

\textit{Step 2.} We use a scaling argument. For every $\t\le s< t\le T$ we show that
$\mathcal{D}_{\frac{1}{\sqrt{t-s}}}\mathcal{A}_{t,s}\mathcal{D}_{\frac{1}{\sqrt{t-s}}}$ coincides
with some matrix $\hat{\mathcal{A}}_{1,0}$ to which we can apply the result of Step 1. We have
\begin{align}
\mathcal{D}_{\frac{1}{\sqrt{t-s}}}\mathcal{A}_{t,s}\mathcal{D}_{\frac{1}{\sqrt{t-s}}}&=
\int_s^t \left(\mathcal{D}_{\frac{1}{\sqrt{t-s}}}\mathcal{E}_{t,\r}\mathcal{D}_{\sqrt{t-s}}\mathbf{e_2}\right)
\left(\mathcal{D}_{\frac{1}{\sqrt{t-s}}}\mathcal{E}_{t,\r}\mathcal{D}_{\sqrt{t-s}}\mathbf{e_2}\right)^{\ast}\frac{d\r}{t-s}\\
&=\int_0^1 \left(\hat{\mathcal{E}}_{1,\r}^{t,s}\mathbf{e_2}\right)\left(\hat{\mathcal{E}}_{1,\r}^{t,s}\mathbf{e_2}\right)^{\ast}
=:\hat{\mathcal{A}}_{1,0}^{t,s}
\end{align}
where
$$\hat{\mathcal{E}}_{\r_1,\r_2}^{t,s}=\mathcal{D}_{\frac{1}{\sqrt{t-s}}}\mathcal{R}_{s+\r_1(t-s),s+\r_1(t-s)}\mathcal{D}_{\sqrt{t-s}},$$
solves the differential system
\begin{align}
 \p_{\r_1}\hat{\mathcal{E}}_{\r_1,\r_2}^{t,s}
 &=(t-s)\mathcal{D}_{\frac{1}{\sqrt{t-s}}}\mathcal{Y}_{s+\r_1(t-s)}\mathcal{D}_{\sqrt{t-s}}\hat{\mathcal{E}}_{\r_1,\r_2}^{t,s}
 =:\hat{\mathcal{Y}}^{t,s}_{\r_1}\hat{\mathcal{E}}_{\r_1,\r_2}^{t,s}
\end{align}
with $\hat{\mathcal{E}}_{\r,\r}^{t,s}=I_2$. A direct computation shows that
  $$(\hat{\mathcal{Y}}^{t,s}_{\r})_{1,2}=(\mathcal{Y}_{s+\r(t-s)})_{1,2}\in [\l^{-1},\l], \qquad
  \|\hat{\mathcal{Y}}^{t,s}_{\r}\|_{\infty}\le (1+T^2)\|\mathcal{Y}_{\r}\|_{\infty}.$$
Therefore \eqref{e34} holds for $\hat{\mathcal{A}}^{t,s}_{1,0}$, uniformly in $t,s$, with $c$
dependent only on $\l$ and $T$.
\endproof

\begin{remark}
Since, for $\t\le s< t\le T$, $\AA_{t,s}^{s,\z}$ is a symmetric and positive definite
matrix, \eqref{e36} also yields an analogous estimate for the inverse: we have
\begin{equation}\label{e36bis}
 c^{-1}\left|\mathcal{D}_{\frac{1}{\sqrt{t-s}}}z\right|^2\le \langle(\AA_{t,s}^{s,\z})^{-1}z,z \rangle\le
 c\left|\mathcal{D}_{\frac{1}{\sqrt{t-s}}}z\right|^2, \qquad \t\le s<t\le T,\ z,\z\in\R^2.
\end{equation}
\end{remark}
The following result is a standard consequence of \eqref{e36} and \eqref{e36bis} (cf., for
instance, Proposition 3.1 in \cite{DiFrancescoPascucci2}).
\begin{proposition}\label{p1}
There exists a positive constant $c$, only dependent on 
$\l_1$, $\l_2$ and $T$, such that
\begin{equation}\label{e_p1}
c^{-1}\G^{\text{\rm heat}}\left(c^{-1} \mathcal{D}_{t-s}, z-\gvar_{t}^{s,\z}\right)\le
Z(t,z;s,\z)\le c\,\G^{\text{\rm heat}}\left(c \mathcal{D}_{t-s}, z-\gvar_{t}^{s,\z}\right),
\end{equation}
for every $\t\le s < t\le T$ and $z,\z\in\R^2$.
\end{proposition}
\begin{remark}{Since $Q_t=\mathcal{D}_{\sqrt{t}}Q_1\mathcal{D}_{\sqrt{t}}$, where $Q_1$ is symmetric and positive definite, estimate \eqref{e_p1} equally holds by replacing $\mathcal{D}_{t-s}$ with $Q_{t-s}$. }
\end{remark}
Next we prove some estimate for the derivatives of $Z(t,z;s,\z)$. We start with the following
\begin{lemma}
 We have
\begin{align}
 (t-s)^{2-i}\left|\left((\AA_{t,s}^{s,\z})^{-1}w\right)_i\right|&\le
 \frac{c}{\sqrt{t-s}}\left|\mathcal{D}_{\frac{1}{\sqrt{t-s}}}w\right|, \label{e37}\\
 (t-s)^{4-i-j}\left|\left((\AA_{t,s}^{s,\z})^{-1}\right)_{ij}\right|&\le
 \frac{c}{t-s}\label{e37bis}
\end{align}
for every $i,j\in\{1,2\}$, $\t\le s <t\le T$ and $w,\z\in\R^2$.
\end{lemma}
\proof We have
\begin{align}
 (t-s)^{2-i}\left|\left((\AA_{t,s}^{s,\z})^{-1}w\right)_i \right|&=\frac{1}{\sqrt{t-s}}
 \left|\left(\mathcal{D}_{\sqrt{t-s}}(\AA_{t,s}^{s,\z})^{-1}\mathcal{D}_{\sqrt{t-s}}\mathcal{D}_{\frac{1}{\sqrt{t-s}}}w\right)_i
 \right|\\ &\le
 \frac{1}{\sqrt{t-s}}\left\|\mathcal{D}_{\sqrt{t-s}}(\AA_{t,s}^{s,\z})^{-1}\mathcal{D}_{\sqrt{t-s}}\right\|\left|\mathcal{D}_{\frac{1}{\sqrt{t-s}}}w\right|.
\end{align}
In order to get \eqref{e37} it suffice to notice that, by \eqref{e36bis}, we have
\begin{align}
\left\|\mathcal{D}_{\sqrt{t-s}}(\AA_{t,s}^{s,\z})^{-1}\mathcal{D}_{\sqrt{t-s}}\right\|\le c.
\end{align}
Taking $w=\mathbf{e}_j$ we also get \eqref{e37bis}.
\endproof

We are ready to state the last result for this section, which is a standard consequence of
estimates \eqref{e37}, \eqref{e37bis} and Proposition \ref{p1} (cf., for instance, Proposition 3.6
in \cite{DiFrancescoPascucci2}).

\begin{proposition}\label{p2}
{There exists a positive constant $c$, only dependent on $\l_1$, $\l_2$ and $T$ such that}
\begin{align}
|\p_x Z(t,z;s,\z)|&\le \frac{{c}}{(t-s)^{\frac{3}{2}}}\G^{\text{\rm heat}}\left({c}
\mathcal{D}_{t-s}, z-\gvar_{t}^{s,\z}\right),\label{e38}\\ |\p_v Z(t,x,v;s,\z)|&\le
\frac{{c}}{\sqrt{t-s}}\G^{\text{\rm heat}}\left({c} \mathcal{D}_{t-s},
z-\gvar_{t}^{s,\z}\right),\label{e38bis}\\ |\p_{vv} Z(t,x,v;s,\z)|&\le
\frac{{c}}{t-s}\G^{\text{\rm heat}}\left({c} \mathcal{D}_{t-s},
z-\gvar_{t}^{s,\z}\right),\label{e38tris}
\end{align}
for every $\t\le s < t\le T$ and $z=(x,v),\z\in\R^2$. 
\end{proposition}

\subsubsection{Upper bound for the fundamental solution}\label{upper}
In this section we assume $\t=0$ for simplicity. We start with some preliminary lemmas. 
\begin{lemma}[Reproduction formula]\label{lemma5} For any $c',c''>0$ we have
\begin{align}
 \Lambda(c',c'')^{-1}\G^{\text{\rm heat}}\left(\frac{c'\wedge c''}{2} \mathcal{D}_{t-s},
 \z''-\z'\right)&\le \int_{\R^2}\G^{\text{\rm heat}}\left(c' \mathcal{D}_{t-\r}, \z'-\y\right)
 \G^{\text{\rm heat}}\left(c'' \mathcal{D}_{\r-s}, \y-\z''\right)d\y\\ &\le
 \Lambda(c',c'')\G^{\text{\rm heat}}\left((c'\vee c'') \mathcal{D}_{t-s}, \z''-\z'\right),
\end{align}
for every $0\le s < \r< t\le T$, $\z',\z''\in\R^2$, where
$\Lambda(c',c'')=\sqrt{\frac{2(c'\vee c'')}{c'\wedge c''}}$. 
\end{lemma}
\proof It is a direct consequence (see also \cite{MR2659772}, Lemma B.1) of the following trivial
estimate
 $$\frac{c'\wedge c''}{2}\mathcal{D}_{t-s}\le c'\mathcal{D}_{t-\r}+c''\mathcal{D}_{\r-s} \le (c'\vee c'')\mathcal{D}_{t-s}.$$
\endproof

\begin{remark}\label{rema1}
Let $\t=0$, $T=1$. If $\hat{Y}$ is a vector field satisfying Assumption \ref{As2} and
$\hat{\g}_t$ is the integral curve
\begin{equation}
    \hat{\g}_t(z)=z+\int_{0}^{t}\hat{Y}_s(\hat{\g}_s(z))ds,\qquad t\in [0,1],
\end{equation}
then $\hat{\g}_1(\cdot)$ is a diffeomorphism of $\R^{2}$. Moreover, since $\hat{Y}$ is Lipschitz
continuous, we have
\begin{equation}\label{e52}
 m^{-1}|z-\hat{\g}_1(\z)|\le |\hat{\g}^{-1}_1(z)-\z|\le m |z-\hat{\g}_1(\z)|, \qquad z,\z\in\R^2,
\end{equation}
for a constant $m$ which depends only on $\l_{2}$.
\end{remark}

\begin{lemma}\label{lemma4}
Let $\gvar_{s}^{t,z}$ be as in \eqref{e31}. There exists a positive constant $m$, only dependent
on $\l_2$ and $T$, such that
 $$m^{-1}\left|\mathcal{D}_{\frac{1}{\sqrt{t-s}}}(z-\gvar_{t}^{s,\z})\right|\le
 \left|\mathcal{D}_{\frac{1}{\sqrt{t-s}}}(\gvar_{s}^{t,z}-\z)\right|\le m
 \left|\mathcal{D}_{\frac{1}{\sqrt{t-s}}}(z-\gvar_{t}^{s,\z})\right|,$$
for every $0\le s< t\le T$ and $z,\z\in\R^2$. 
\end{lemma}
\proof We use again a scaling argument: we set $z'=\mathcal{D}_{\sqrt{t-s}}z$ and
 $$\hat{\gvar}_{\r}(z)=\mathcal{D}_{\frac{1}{\sqrt{t-s}}}\gvar_{s+\r (t-s)}^{s,z'},\quad
 \hat{Y}_{\r}(z)=(t-s)\mathcal{D}_{\frac{1}{\sqrt{t-s}}}Y_{s+\r(t-s)}(z'),\qquad  \r\in [0,1].$$
Then we have
 $$\hat{\gvar}_{\r}(z)=z+\int_{0}^{\r}\hat{Y}_{u}(\hat{\gvar}_{u}(z))du,\qquad  \r\in [0,1].$$
As in the proof of Proposition \ref{p0}, we have that $\hat{Y}$ satisfies Assumption \ref{As2}. By
Remark \ref{rema1}, estimate \eqref{e52} holds for $\hat{\gvar}_{\r}(z)$. {To conclude, it
suffices to substitute $z$ and $\z$ with $\bar{z}=\mathcal{D}_{\frac{1}{\sqrt{t-s}}}z$ and
$\bar{\z}=\mathcal{D}_{\frac{1}{\sqrt{t-s}}}\z$ in \eqref{e52}.}
\endproof
\begin{lemma}\label{p3}
Let $(\mathcal{K}_{t}Z)_{k}$ be as in \eqref{expansion2}. {There exists a constant $c>0$, only dependent on 
$\l_1$, $\l_2$, $\a$ and $T$ such that}
\begin{align}
 \left|(\mathcal{K}_{t}Z)_{k}(t,z;s,\z)\right|\le \frac{M_k}{(t-s)^{1-\frac{k\a}{2}}} \G^{\text{\rm
 heat}}\left({c} \, m^k \mathcal{D}_{t-s}, z-\gvar_{t}^{s,\z}\right), \qquad 0\le s < t\le T, \
 z,\z\in\R^2,
\end{align}
where $m$ is the constant in Lemma \ref{lemma4} and $M_k=2^{\frac{k-1}{2}}{c}^k
m^{q_k}\frac{\G^k_E(\frac{\a}{2})}{\G_E(\frac{k\a}{2})}$, with $q_1=0, \ q_2=\frac{1}{2},
q_k=q_{k-1}+\frac{k-2}{2}$ for $k\ge 2$.
\end{lemma}
\proof We give the proof for $k=1$. 
{The general case follows by induction, exploiting Lemmas \ref{lemma4} and \ref{lemma5} as in the
proof of estimate \eqref{u_bound}.}
\begin{align}
 (\mathcal{K}_{t}Z)_1(t,z;s,\z)&= (\mathcal{L}_{t}-\mathcal{L}_{t}^{s,\z})Z(t,z;s,\z)\\
 &=\frac{1}{2}\left(a_{t}(z)-a_{t}(\gvar_{t}^{s,\z})\right)\p_{vv}Z(t,z;s,\z)
 {+b_{t}(z)\p_{v}Z(t,z;s,\z)}+\\
 &\quad +\langle
 Y_{t}(z)-\bar{Y}_{t}^{s,\z}(z),\nabla Z(t,z;s,\z)\rangle\\
 &=:E_1+E_2+E_3.
\end{align}
By Assumption {\ref{As1}} and Proposition \ref{p2} we have
\begin{align}
|E_1|&\le \frac{{c}}{t-s} |z-\gvar_{t}^{s,\z}|^{\a}\G^{\text{\rm heat}}\left({c}
\mathcal{D}_{t-s}, z-\gvar_{t}^{s,\z}\right)\\ &\le \frac{{c'}
}{(t-s)^{1-\frac{\a}{2}}}\left|\mathcal{D}_{\frac{1}{\sqrt{t-s}}}(z-\gvar_{t}^{s,\z})\right|^{\a}
\G^{\text{\rm heat}}\left({c} \mathcal{D}_{t-s}, z-\gvar_{t}^{s,\z}\right)\le\\ \intertext{(by
\eqref{e43})} &\le \frac{{c''}
}{(t-s)^{1-\frac{\a}{2}}}\G^{\text{\rm heat}}\left({c''}\mathcal{D}_{t-s},
z-\gvar_{t}^{s,\z}\right).
\end{align}
{{By Assumption {\ref{As1}} and Proposition \ref{p2} we also have $$|E_2|\le
\frac{{c}}{\sqrt{t-s}}\G^{\text{\rm heat}}\left({c} \mathcal{D}_{t-s},
z-\gvar_{t}^{s,\z}\right).$$}} As for $E_3$, we have
\begin{align}
 |(Y_{t}(z)-\bar{Y}_{t}^{s,\z}(z))_1|=|Y_{1,t}(z)-Y_{1,t}(\gvar_{t}^{s,\z})-\p_vY_{1,t}(\gvar_{t}^{s,\z})(z-\gvar_{t}^{s,\z})|
 \le c |z-\gvar_{t}^{s,\z}|^{1+\a},
\end{align}
because $\p_vY_{1,t}$ is H\"older continuous by Assumption \ref{As2}: here we use the elementary
inequality
 $$\left|\int_0^1 (f'(y+t(x-y))-f'(y))(x-y)dt \right|\le c_{{\a}}|x-y|^{1+\a}.$$
which is valid for $f\in C^{1+\a}$. 
On the other hand, we have
\begin{align}
 |(Y_{t}(z)-\bar{Y}_{t}^{s,\z}(z))_2|\le c |z-\gvar_{t}^{s,\z}|.
\end{align}
Therefore, by Proposition \ref{p2}, we have
\begin{align}
|E_3|&\le
{c}\left(\frac{1}{(t-s)^{\frac{3}{2}}}|z-\gvar_{t}^{s,\z}|^{1+\a}+\frac{1}{(t-s)^{\frac{1}{2}}}|z-\gvar_{t}^{s,\z}|\right)
\G^{\text{\rm heat}}\left({c}\mathcal{D}_{t-s}, z-\gvar_{t}^{s,\z}\right) \\ &\le
\frac{{c'}}{(t-s)^{1-\frac{\a}{2}}}
\left(\left|\mathcal{D}_{\frac{1}{\sqrt{t-s}}}(z-\gvar_{t}^{s,\z})\right|^{1+\a}+
\left|\mathcal{D}_{\frac{1}{\sqrt{t-s}}}(z-\gvar_{t}^{s,\z})\right|\right) \G^{\text{\rm
heat}}\left({c}\mathcal{D}_{t-s}, z-\gvar_{t}^{s,\z}\right)\le \intertext{(by \eqref{e43})} &\le
\frac{{c''}}{(t-s)^{1-\frac{\a}{2}}} \G^{\text{\rm heat}}\left({c''}\mathcal{D}_{t-s},
z-\gvar_{t}^{s,\z}\right).
\end{align}
\endproof
The following result is proved in \cite{MR2659772}, Proposition 5.2.
\begin{lemma}\label{p4} For any $\e>0$ there exist a positive constant $c$, only dependent on $\l_{1},\l_{2},\a,T$ and $\e$, such that
\begin{equation}
 \int_{\R^2}\pG(t,z;\r,\y)(\r-s)^2\G^{\text{\rm heat}}(\e \mathcal{D}_{\r-s},\y-\gvar_{\r}^{s,\z})d\y\le
 c\, \G^{\text{\rm heat}}\left(c \mathcal{D}_{t-s}, z-\gvar_{t}^{s,\z}\right),
\end{equation}
for $s<\r<t\le T$ and $z,\z\in \R^2$.
\end{lemma}


We close this section by proving the Gaussian upper bound in \eqref{e40}. Consider the parametrix
expansion \eqref{expansion} with $0<t-s\le 1$. By Proposition \ref{p1}, the first term in the RHS
of \eqref{expansion} is bounded by $c\,\G^{\text{\rm heat}}(c\mathcal{D}_{t-s},
z-\g_t^{s,\z})$. On the other hand, if $N\ge\frac{6}{\a}$ 
then $(\r-s)^{1-\frac{N\a}{2}}\le (\r-s)^2$ and therefore the last term in the RHS of
\eqref{expansion} is bounded by the same quantity, by Lemmas \ref{p3} and \ref{p4}.

Finally, 
denoting with $c_k$ a positive constant dependent on $\l_{1},\l_{2},\a,T$ and $k$, we have
\begin{align}
 \int_{s}^t\int_{\R^2}&Z(t,z;\r,\y)(\mathcal{K}_{\r}Z)_k(\r,\y;s,\z)d\y d\r\le
\intertext{(by Lemmas \ref{p3} and \ref{p4})}
 & \quad \le c_k
 \int_{s}^t (\r-s)^{\frac{k\a}{2}-1}\int_{\R^2}\G^{\text{\rm heat}}(c\mathcal{D}_{t-\r},
 z-\g_t^{\r,\y})\G^{\text{\rm heat}}(c_k \mathcal{D}_{\r-s}, \y-\g_{\r}^{s,\z})d\y d\r \le
 \intertext{(by Lemma \ref{lemma4})} & \quad \le c_k \int_{s}^t
 \r^{\frac{k\a}{2}-1}\int_{\R^2}\G^{\text{\rm heat}}(c'\mathcal{D}_{t-\r},
 \gvar_{\r}^{t,z}-\y)\G^{\text{\rm heat}}(c_k \mathcal{D}_{\r-s}, \y-\g_{\r}^{s,\z})d\y d\r  \le
 \intertext{(by Lemma \ref{lemma5})} & \quad \le c_k \int_{s}^t \r^{\frac{k\a}{2}-1}\G^{\text{\rm
 heat}}(c_k \mathcal{D}_{t-s},\gvar_{\r}^{t,z}-\g_{\r}^{s,\z})d\r \le \intertext{(again by Lemma
 \ref{lemma4})} & \quad \le c_k \G^{\text{\rm heat}}(c_k \mathcal{D}_{1},z-\g_{t}^{s,\z}).\label{u_bound}
\end{align}
This proves the upper bound for $0<t-s\le 1$. The general case can be recovered by a scaling
argument, similar to that of Proposition \ref{p0}.

\subsubsection{Lower bound for the fundamental solution}
We first derive a local bound, starting from the parametrix expansion \eqref{expansion0} and
exploiting the results of Section \ref{upper}. We have
\begin{align}
 \pG(t,z;s,\z)&\ge
 Z(t,z;s,\z)-\int_s^t\int_{\R^2}\pG(t,z;\r,\y)\left|\mathcal{K}_{\r}Z(\r,\y;s,\z)\right|d\y
 d\r\ge
\intertext{(by Lemmas \ref{p1} and \ref{p3} and the upper bound \eqref{e40})}
 &\ge c^{-1}\G^{\text{\rm heat}}(c^{-1}\mathcal{D}_{t-s},z-\gvar_t^{s,\z})\\
 &\quad -\int_s^t\frac{c}{(\r-s)^{1-\frac{\a}{2}}} \int_{\R^2}\G^{\text{\rm
 heat}}(c\mathcal{D}_{t-\r},z-\gvar_t^{\r,\y})\G^{\text{\rm
 heat}}(c\mathcal{D}_{\r-s},z-\gvar_{\r}^{s,\z})d\y d\r\ge
\intertext{(by Lemma \ref{lemma5})}
 &\ge  c^{-1}\G^{\text{\rm
 heat}}(c^{-1}\mathcal{D}_{t-s},z-\gvar_t^{s,\z})-\frac{c}{2}(t-s)^{\frac{\a}{2}}\G^{\text{\rm
 heat}}(c\mathcal{D}_{t-s},z-\gvar_t^{s,\z}).
\end{align}
Let $d_{t_2,t_1}(z_2,z_1):=|\mathcal{D}_{t_2-t_1}(z_2-\gvar_{t_2}^{t_1,z_1})|$ denote the
``control metric'' of the system. A direct computation shows that $\G^{\text{\rm
heat}}(c\mathcal{D}_{t-s},z-\gvar_t^{s,\z})\le \G^{\text{\rm
heat}}(c^{-1}\mathcal{D}_{t-s},z-\gvar_t^{s,\z})$ if $d_{t,s}(z,\z)\le \varrho_c$ where
$\varrho_c=\sqrt{\frac{2c\ln c}{c^2-1}}$. Then we have
\begin{align}\label{e60}
 \pG(t,z;s,\z)&\ge \left(\frac{1}{c^{2}}-\frac{(t-s)^{\frac{\a}{2}}}{2}\right)\G^{\text{\rm
 heat}}(c^{-1}\mathcal{D}_{t-s},z-\gvar_t^{s,\z})
 \ge \frac{1}{2c}\G^{\text{\rm
 heat}}(c^{-1}\mathcal{D}_{t-s},z-\gvar_t^{s,\z})
\end{align}
if $d_{t,s}(z,\z)\le \varrho_c$ and $0< t-s\le T_c:=c^{-\frac{4}{\a}}$.

In order to pass from the local to the global bound, we use a chaining procedure: we first need to
define a sequence of points $(t_k, z_k)$ such that $t_0=s, z_0=\z, t_{M+1}=t, z_{M+1}=z$ for some
integer $M$ (to be defined later), along which we can control the increments with respect to the
control metric $d_{t_{k-1},t_k}(z_{k+1},z_k)$. Let us consider the controlled version of the
system \eqref{e31}:
 $$\psi_{\r}^{s,\z}=\z+\int_s^{\r}\left(Y_{\theta}(\psi_{\theta}^{s,\z})+v_s\mathbf{e}_2\right)d\theta,
 \qquad \r\in [s,t].$$
We have the following (see \cite{Polidoro1}, \cite{MR2231876} and \cite{MR2659772}, Propositions
4.1 and 4.2):
\begin{lemma}\label{lemma6} There exists a control $(v_{\r})_{s\le \r\le t}$ with values in $\R^2$ such that
\begin{itemize}
\item[i)] the solution $\psi_{\r}^{s,\z}$ associated with $v_{\r}$ reaches $z$ at time $t$, that is $\psi_t^{s,\z}=z$;
\item[ii)] there exist two constants $m_1,m_2> 0$, only dependent on the constants of Assumptions \ref{As1}-\ref{As2}, such that
\begin{align}
  \int_s^t|v_{\r}|^2d\r \ge m_1 \left|\mathcal{D}_{\frac{1}{\sqrt{t-s}}}(z-\gvar_t^{s,\z})\right|^2,\qquad
  \sup_{s\le \r\le t}|v_{\r}|^2\le \frac{m_2}{t-s}\left|\mathcal{D}_{\frac{1}{\sqrt{t-s}}}(z-\gvar_t^{s,\z})\right|^2.
\end{align}
\end{itemize}
\end{lemma}
We set $$t_i=s+i\frac{t-s}{M+1}=s+i\e, \quad z_k=\psi_{r_k}^{s,\z}, \qquad i=1,\cdots,M,$$ where
$\psi_{\r}^{s,\z}$ is the optimal path of Lemma \ref{lemma6} and $M$ is the smallest integer
greater than
\begin{equation}
\max\left\{\frac{K^2 d^2_{t,s}(z,\z)}{\varrho_c^2}, \frac{T}{T_c}\right\}.
\end{equation}
with $K=\frac{12 m^2 m_2}{m_1}$, where $m$, $m_1$ and $m_2$ are the constants in Lemmas \ref{lemma4} and \ref{lemma6}.
Finally we define the sets
 $$B_i(r):=\{z\in\R^2 \mid |\mathcal{D}_{\frac{1}{\sqrt{\e}}}(z-\gvar_{t_i}^{t_{i-1},z_{i-1}})|
 +|\mathcal{D}_{\frac{1}{\sqrt{\e}}}(z_{i+1}-\gvar_{t+i}^{t_{i},z})|\le r\},$$
and write
\begin{equation}\label{e61}
\pG(t,z;s,\z)\ge \int_{B_1(\varrho_c/3)}\!\!\cdots\int_{B_M(\varrho_c/3)}
\pG(t,z;t_M,\z_M)\prod_{j=1}^{M-1}\pG(t_{j+1},\z_{j+1};t_j,\z_j)\pG(t_{1},\z_{1};s,\z)d\z_1\dots
d\z_M.
\end{equation}
By definition of $M$ we have
  $$t_{j+1}-t_j=\frac{t-s}{M+1}\le \frac{T}{M+1}\le T_c.$$
On the other hand, if $\z_i\in
B_i\left(\frac{\varrho_c}{3}\right)$ for $i=1,\dots , M-1$ we have
\begin{align}
 d_{t_{i+1},t_i}(\z_{i+1},\z_i)&=|\mathcal{D}_{\frac{1}{\sqrt{\e}}}(\z_{i+1}-\gvar_{t_{i+1}}^{t_{i},\z_{i}})|\\
 &=|\mathcal{D}_{\frac{1}{\sqrt{\e}}}(\z_{i+1}-\gvar_{t_{i+1}}^{t_{i},z_{i}})|+
 |\mathcal{D}_{\frac{1}{\sqrt{\e}}}(z_{i+1}-\gvar_{t_{i+1}}^{t_{i},z_{i}})|+
 |\mathcal{D}_{\frac{1}{\sqrt{\e}}}(z_{i+1}-\gvar_{t_{i+1}}^{t_{i},\z_{i}})|=:E_1+E_2+E_3,
\end{align}
where $E_1+E_3\le \frac{2}{3}\varrho_c$. By Lemma \ref{lemma6}, we have
\begin{align}\label{e64}
 E_2&\le m_1^{-1}\left(\int_{t_i}^{t_{i+1}}|v_{\r}|^2d\r\right)^{\frac{1}{2}}
 \le \frac{m_2}{m_1}\sqrt{\frac{\e}{t-s}}|\mathcal{D}_{\frac{1}{\sqrt{t-s}}}(z-\gvar_t^{s,\z})|
 =\frac{m_2}{m_1}\frac{d_{t,s}(z,\z)}{\sqrt{M+1}}\le \frac{\varrho_c}{12 m^2}.
\end{align}
Therefore $d_{t_{i+1},t_i}(\z_{i+1},\z_i)\le \varrho_c$ and we can use \eqref{e60} repeatedly in
\eqref{e61} to get
\begin{equation}\label{e62}
 \pG(t,z;s,\z)\ge (2c)^{-(M+1)}\left|\prod_{i=1}^M
 B_i\left(\frac{\varrho_c}{3}\right)\right|\left(\frac{c(M+1)^2}{(t-s)^2}\right)^{M+1}
 \exp\left(-\frac{c}{2}\varrho_c^2(M+1)\right).
\end{equation}
Assume for a moment the validity of the inequality
\begin{equation}\label{e63}
 \left|B_i\left(\frac{\varrho_c}{3}\right)\right|\ge C_0 \pi\left(\frac{t-s}{M+1}\right)^2\varrho_c^2
\end{equation}
for some positive constant $C_0$ (only dependent on the constants of Assumptions
\ref{As1}-\ref{As2}). Then we have
\begin{align}
 \pG(t,z;s,\z)\ge C_1 C_2^M\frac{1}{2\pi \sqrt{\det\mathcal{D}_{t-s}}}\exp\left(-\frac{c}{2}\varrho_c^2 M\right)
 \ge \frac{C_3}{2\pi \sqrt{\det\mathcal{D}_{t-s}}}\exp\left(-\frac{C_4}{2}M\right),
\end{align}
for some positive constants $C_1,\dots,C_4$. Now, if $TT^{-1}_c \le \frac{K^2
d^2_{t,s}(z,\z)}{\varrho_c^2}$ and $M< 2 \frac{K^2d^2_{t,s}(z,\z)}{\varrho_c^2}$, we have
 $$\pG(t,z;s,\z)\ge \frac{C_3}{2\pi \sqrt{\det\mathcal{D}_{t-s}}}\exp\left(-\frac{C_5}{2}d^2_{t,s}(z,\z)\right) =C_6 \G^{\text{\rm
 heat}}(C_5^{-1}\mathcal{D}_{t-s},z-\gvar_t^{s,\z}).$$
On the other hand, if $M<2TT_c^{-1}$ then
 $$\pG(t,z;s,\z)\ge \frac{C_7}{2\pi \sqrt{\det\mathcal{D}_{t-s}}}\ge \frac{C_7}{2\pi
 \sqrt{\det\mathcal{D}_{t-s}}}\exp\left(-\frac{C_5}{2}d^2_{t,s}(z,\z)\right)= C_8 \G^{\text{\rm
 heat}}(C_5^{-1}\mathcal{D}_{t-s},z-\gvar_t^{s,\z}),$$
and this proves the lower bound.

We are left with the proof of \eqref{e63}.
Let $\tilde{B}_i(r)=\{z, \ |\mathcal{D}_{\frac{1}{\sqrt{\e}}}(z-z_i)|\le r \}$: a direct
computation shows $|\tilde{B}_i(r)|=\pi\e^2r^2$. Then it is enough to show that
$B_i\left(\frac{\varrho_c}{3}\right)\supseteq \tilde{B}_i\left(\frac{\varrho_c}{C}\right)$ for a
positive constant $C$ (only dependent on $\l_{1},\l_{2},\a$ and $T$). For any $z\in
\tilde{B}_i(r)$ we have
\begin{align}
 &|\mathcal{D}_{\frac{1}{\sqrt{\e}}}(z-\gvar_{t_i}^{t_{i-1},z_{i-1}})|
 +|\mathcal{D}_{\frac{1}{\sqrt{\e}}}(z_{i+1}-\gvar_{t+i}^{t_{i},z})|\le \\
 &\qquad \le |\mathcal{D}_{\frac{1}{\sqrt{\e}}}(z-z_i)|+|\mathcal{D}_{\frac{1}{\sqrt{\e}}}(z_i-\gvar_{t_i}^{t_{i-1},z_{i-1}})|+
 m|\mathcal{D}_{\frac{1}{\sqrt{\e}}}(z-z_i)|+m^2|\mathcal{D}_{\frac{1}{\sqrt{\e}}}(z_{i+1}-\gvar_{t_{i+1}}^{t_{i1},z_{i1}})|\le
\intertext{(by \eqref{e64})}
 &\qquad \le (1+m)r+\frac{\varrho_c}{6}.
\end{align}
Then it is sufficient to take $r\le \frac{\varrho_c}{6(1+m)}$ and this concludes the proof.
\endproof

\subsubsection{Gaussian bounds for $\p_v\pG$ and $\p_{vv}\pG$}
The following lemma provides an alternative representation formula for $\pG$ which will be used to
prove the bounds for the derivatives. As a general rule, until the end of the section we will
always denote with $c$ a positive constant, only dependent on $\l_1, \l_2, \a$ and $T$ in
Assumptions \ref{As1}-\ref{As2}.
\begin{lemma} We have
\begin{equation}
 \pG(t,z;s,\z)=Z(t,z;s,\z)+\int_s^t\int_{\R^2}Z(t,z;\r,\y)\phi(r,\y;s,\z)d\y d\r, \qquad
 \t\le s< t\le T,\ z,\z\in\R^2,
\end{equation}
where
\begin{equation}
 \phi(\cdot,\cdot;s,\z)=\sum_{k\ge 1}(\mathcal{K}Z)_k(\cdot,\cdot;s,\z)
\end{equation}
is uniformly convergent in $(s,T)\times \R^2$. Moreover, there exists a positive constant $c$ such
that
\begin{align}\label{ee50}
 |\phi(t,z;s,\z)|&\le \frac{c}{(t-s)^{1-\frac{\a}{2}}}\pG^{\text{\rm heat}}(c\mathcal{D}_{t-s},z-\gvar_t^{s,\z}),\\
 |\phi(t,z;s,\z)-\phi(t,z';s,\z)|&\le c\frac{{d_{\mathcal{L}}\left((t,z),(t,z')\right)^{\frac{\a}{2}}}}{(t-s)^{1-\frac{\a}{2}}}
 \left(\pG^{\text{\rm heat}}(c\mathcal{D}_{t-s},z-\gvar_t^{s,\z})+\pG^{\text{\rm
 heat}}(c\mathcal{D}_{t-s},z'-\gvar_t^{s,\z})\right), \label{e50bis}
\end{align}
for every $\t\le s< t\le T$ and $z,z',\z\in \R^2$, {where $d_{\mathcal{L}}$ is the intrinsic
distance in \eqref{ee37bis}.}
\end{lemma}
\proof We start from the parametrix representation \eqref{expansion} and show that the remainder
 $$R_N(t,z;s,\z):=\int_s^t\int_{\R^2}\pG(t,z;\r,\y)(\mathcal{K}_{\r}Z)_N(\r,\y;s,\z)d\y d\r$$
converges uniformly to $0$ as $N$ tends to infinity. By the Gaussian upper bound \eqref{e40},
Lemmas \ref{p3} and \ref{lemma5}, we have
\begin{align}
 |R_N(t,z;s,\z)|&\le c M_N\int_s^t\frac{1}{(t-\r)^{1-\frac{N\a}{2}}}\int_{\R^2} \G^{\text{\rm
 heat}}(c\mathcal{D}_{t-\r},\gvar_{\r}^{t,z}-\y)\G^{\text{\rm
 heat}}(c_N\mathcal{D}_{\r-s},\y-\gvar_{s}^{s,z})d\y d\r \\ &\le c
 M_N\int_s^t\frac{1}{(t-\r)^{1-\frac{N\a}{2}}} \G^{\text{\rm
 heat}}(c\mathcal{D}_{t-s},\gvar_{\r}^{t,z}-\gvar_{s}^{s,\z})d\r \\ &\le c
 M_N(t-s)^{-2}\int_s^t\frac{1}{(t-\r)^{1-\frac{N\a}{2}}}d\r \\ &\le c
 \frac{M_N}{N}(t-s)^{\frac{N\a}{2}-2},
\end{align}
with $M_N=c^N m^{q_N}\frac{\G^N_E(\frac{\a}{2})}{\G_E(\frac{N\a}{2})}$, converges to zero by the
properties of the Euler Gamma function $\G_E$.

Next, exploiting the lower bound for $\pG$ we can replace the Gaussian function $\G^{\text{\rm
heat}}$ in Propositions \ref{p1} and \ref{p2} by an appropriate fundamental solution
satisfying an exact reproduction formula. Then, repeating the arguments in the proof of Lemma
\ref{p3}, we get
\begin{equation}
 \left|(\mathcal{K}_{t}Z)_{k}(t,z;s,\z)\right|\le \frac{M_k}{(t-s)^{1-\frac{k\a}{2}}} \G^{\text{\rm
 heat}}\left(c \mathcal{D}_{t-s}, z-\gvar_{t}^{s,\z}\right), \qquad \t\le s < t\le T, \
 z,\z\in\R^2.
\end{equation}
Then estimate \eqref{ee50} easily follows. Estimate \eqref{e50bis} can be proved by standard
arguments (see, for instance, Lemma 6.1 in \cite{DiFrancescoPascucci2}).
\endproof

Now we show that
\begin{align}
\left|\int_s^t\int_{\R^2}\p_vZ(t,z;\r,\y)\phi(r,\y;s,\z)d\y d\r\right| &\le
\frac{c}{\sqrt{t-s}}\G^{\text{\rm heat}}\left(c \mathcal{D}_{t-s},z-\gvar_{t}^{s,\z}\right),
\label{e51}\\ \left|\int_s^t\int_{\R^2}\p_{vv}Z(t,z;\r,\y)\phi(r,\y;s,\z)d\y d\r\right| &\le
\frac{c}{t-s}\G^{\text{\rm heat}}\left(c  \mathcal{D}_{t-s},z-\gvar_{t}^{s,\z}\right),
\label{e51bis}
\end{align}
for $\t\le s <t\le T$ and $z,\z\in\R^2$. Formula \eqref{e51} is a standard consequence of Lemma
\ref{lemma5} and estimates \eqref{e38bis} and \eqref{e50}. Estimate \eqref{e51} is less obvious.
We have
\begin{align}
 \int_{\R^2}\p_{vv}Z(t,z;\r,\y)\phi(\r,\y;s,\z)s\y&=
 \int_{\R^2}\p_{vv}Z(t,z;\r,\y)(\phi(\r,\y;s,\z)-\phi(\r,w;s,\z))d\y\\
 &\quad+\phi(\r,w;s,\z)\int_{\R^2}\p_{vv}(Z(t,z;\r,\y)-\pG_{\r,w}(t,z;\r,\y))d\y\\
 &\quad+\phi(\r,w;s,\z)\int_{\R^2}\p_{vv}\pG_{\r,w}(t,z;\r,\y)d\y\\
 &=:I_1+I_2+I_3.
\end{align}


Then, by choosing $w=\gvar_{\r}^{t,z}$ we can rely on the H\"older regularity of $\phi$ and
$\pG_{\r,y}$ to remove the singularity in $t=\r$. Here we show how to handle $I_1$ in detail: by
estimates \eqref{e38tris} and \eqref{e50bis} we have
\begin{align}
 |I_1|\le & \frac{c}{(\r-s)^{1-\frac{\a}{2}}} \int_{\R^2}\frac{{d_{\mathcal{L}}\left((\r,\gvar_{\r}^{t,z}),(\r,\y)\right)^{\a}}}{t-\r}
 \pG^{\text{\rm heat}}(c \mathcal{D}_{t-\r},z-\gvar_t^{\r,\y})\times \\
 &\times\underbrace{\left(\pG^{\text{\rm heat}}(c
 \mathcal{D}_{\r-s},\y-\gvar_{\r}^{s,\z})+\pG^{\text{\rm heat}}(c
 \mathcal{D}_{\r-s},\gvar_{\r}^{t,z}-\gvar_{\r}^{s,\z})\right)}_{=:J(\y)}d\y \\
 \le & \frac{c}{(t-\r)^{1-\frac{\a}{2}}(\r-s)^{1-\frac{\a}{2}}} \int_{\R^2}
 {\left|\left(0,\mathcal{D}_{\frac{1}{\sqrt{t-\r}}}(z-\gvar_{t}^{\r,\y})\right)\right|^{\a}_{\mathcal{L}}}\pG^{\text{\rm
 heat}}(c \mathcal{D}_{t-\r},z-\gvar_t^{\r,\y}) J(\y)d\y\le
\intertext{(by \eqref{e43})}
  \le &
 \frac{c}{(t-\r)^{1-\frac{\a}{2}}(\r-s)^{1-\frac{\a}{2}}} \int_{\R^2}\pG^{\text{\rm heat}}(c
 \mathcal{D}_{t-\r},z-\gvar_t^{\r,\y})J(\y)d\y \\ =
 &\frac{c}{(t-\r)^{1-\frac{\a}{2}}(\r-s)^{1-\frac{\a}{2}}}\left(I_{11}+I_{12}\right)
\end{align}
where $$I_{11}\le c\,\pG^{\text{\rm heat}}(c \mathcal{D}_{t-s},z-\gvar_t^{s,\z}),$$ by Lemma
\ref{lemma5}, and
\begin{align}
 I_{12}\le \pG^{\text{\rm heat}}(c \mathcal{D}_{\r-s},z-\gvar_t^{s,\z})\int_{\R^2}\pG^{\text{\rm
 heat}}(c \mathcal{D}_{t-\r},z-\gvar_t^{\r,\y})d\y \le c\, \pG^{\text{\rm heat}}(c
 \mathcal{D}_{t-s},z-\gvar_t^{s,\z}),
\end{align}
because the integral is bounded by a constant and the matrix $\mathcal{D}_{\r-s}$ is
increasing in $\r$. $I_2$ can be treated similarly, once we notice that
\begin{equation}
\left|\p_{vv}\pG_{s,y}(t,z;s,\z)-\p_{vv}\pG_{s,{w}}(t,z;s,\z)\right|\le
c\frac{{d_{\mathcal{L}}\left((s,y),(s,w)\right)^{\a}}}{t-s}\pG^{\text{\rm heat}}(c
\mathcal{D}_{t-s},z-\gvar_t^{s,\z}),
\end{equation}
for $\t\le s <t\le T$ and $z,\z,y,w\in\R^2$ (see also \cite{DiFrancescoPascucci2}, Lemma 5.2).
{Lastly, $I_3=0$: indeed, for every $s< \r < t$ and $w\in\R^2$ we have
$\int_{\R^2}\pG_{\r,w}(t,z;\r,\y)d\y=1$ and therefore
$$\p_{vv}\int_{\R^2}\pG_{\r,w}(t,z;\r,\y)d\y=0.$$} Integrating in $\r$ over the interval $(s,t)$
we get estimate \eqref{e51bis}.

\section{Finale: proof of Theorem \ref{t1}}\label{proof}
For any fixed $\t\in [0,T)$ and {$\o\in\O$}, let $\mathcal{K}_{\t}$ the operator of the form
\eqref{pde0}, as defined by \eqref{PDE} and \eqref{e04} through the random change of variable
$\gIW_{\t,t}$. By Assumptions \ref{Ass1}-\ref{Ass2} and Lemma \ref{lemma1}, $\mathcal{K}_{\t}$
satisfies Assumptions \ref{As1}-\ref{As2} for a.e. $\w\in\O$. Then, by Theorem \ref{t3},
$\mathcal{K}_{\t}$ admits a fundamental solution $\pG_{\t}$: we set
\begin{equation}\label{ee43}
 \mathbf{\pG}(t,x,v;\t,\z)=\pG_{\t}(t,x,\gIWi_{t,\t}(x,v);\t,\z), \qquad \t<t\le T, \ x,v \in \R, \ z\in \R^2.
\end{equation}
Combining Theorems \ref{th1}, \ref{t3} and Lemma \ref{lemma1} we infer that
$\mathbf{\pG}(\cdot,\cdot,\cdot;\t,\z)\in \mathbf{C}_{t_0,T}^0$ for any $t_0\in (\t,T]$, is twice
continuously differentiable in the variable $v$ and satisfies \eqref{e50} with probability one.
Now, for any bounded and continuous function $\phi$ and $z_0\in\R^2$, we have
\begin{align}
 \int_{\R^2}\mathbf{\pG}(t,z;\t,\z)\phi(\z)d\z-\phi(z_0)&=
  \int_{\R^2}\pG_{\t}(t,z;\t,\z)\phi(\z)d\z-\phi(z_0)+\\
 &\quad+\int_{\R^2}\left(\pG_{\t}(t,x,\gIWi_{t,\t}(x,v);\t,\z)-\pG_{\t}(t,z;\t,\z)\right)\phi(\z)d\z= \\
 &= I_1(t,z,\t)+I_2(t,z,\t).
\end{align}
Now, by Theorem \ref{t3} and the dominated convergence theorem, we have
 $$\lim_{(t,z)\to(\t,z_{0})\atop t>\t}I_i(t,z,\t)=0,\qquad i=1,2.$$
This proves the first part of the thesis.

 The Gaussian bounds \eqref{e31bbis} follow directly from the definition \eqref{ee43} and the analogous
estimates \eqref{e40} for $\pG_{\t}$ in Theorem \ref{t3}. Moreover, since
 $$\p_{v}\mathbf{\pG}(t,x;\t,\x)=(\p_v \pG_{\t})\left(t,x,\gIWi_{t,\t}(x,v);
  \t,\z\right)\p_v\gIWi_{t,\t}(x,v),$$
the gradient estimate \eqref{e32bbis} follows from the analogous estimate \eqref{e41} for
$\pG_{\t}$ and from Lemma \ref{lemma1}. The proof of \eqref{e33bbis} is analogous.
\endproof

\bibliographystyle{imsart-nameyear}
\bibliography{bib}
\end{document}